\documentclass[11pt,reqno]{amsart}
\usepackage{fixltx2e}                    %% Get the latest LaTeX enhancements
\usepackage{mathpazo}  % This provides the palatino font and some nice "old style" numbers.
\usepackage[utf8]{inputenc}            %% Customize to your encoding !!!
\usepackage{amsmath}                     %% Imprescindible
\usepackage{amssymb, latexsym, stmaryrd, amsthm, dsfont, amsfonts, amsbsy,amsthm, amsmath, mathrsfs}            %% Imprescindible
\usepackage{mathtools}                   %% for some math typesetting tricks
\usepackage{bm}                          %% for good bold math symbols
\usepackage{enumerate}                   %% enhanced "enumerate" environments
\usepackage{verbatim}                    %% for enhanced "comment" environment
\usepackage{url}   
\usepackage{lscape}                      %% for \url command in bibliography
\usepackage{soul}
\usepackage{tikz} 

\usepackage[margin=1.3in]{geometry}

\usepackage{microtype}  
\usepackage[all, knot]{xy}
        \xyoption{arc} 
        \xyoption{web}                 %% for the best typographic result
\makeatletter                            %% A correction needed for microtype:
\def\MT@register@subst@font{\MT@exp@one@n\MT@in@clist\font@name\MT@font@list
 \ifMT@inlist@\else\xdef\MT@font@list{\MT@font@list\font@name,}\fi}
\makeatother

\usepackage[pdftex,bookmarks,bookmarksnumbered,linktocpage,   %%  customize to your liking
         colorlinks,linkcolor=blue,citecolor=blue]{hyperref}

\newcommand{\bit}{\begin{itemize}}    % but see also \benbullet below
\newcommand{\eit}{\end{itemize}}
\newcommand{\ben}{\begin{enumerate}}
\newcommand{\een}{\end{enumerate}}
\newcommand{\benormal}{\ben[\normalfont 1.]}   % *
\let\enormal\een
\newcommand{\benroman}{\ben[\normalfont (i)]}  % *
\let\eroman\een
\newcommand{\bde}{\begin{description}}
\newcommand{\ede}{\end{description}}

\newcommand{\?}{\ensuremath{\mkern0.4\thinmuskip}}   % very small math space
\let\leq=\leqslant
\let\nleq=\nleqslant
\let\geq=\geqslant
\let\epsilon=\varepsilon
\let\Lambda\varLambda
\let\Gamma\varGamma
\let\Delta\varDelta
\let\Lambda\varLambda
\let\omega\varOmega
\let\Theta\varTheta
\let\Xi\varXi
\let\Pi\varPi
\let\Sigma\varSigma
\let\class=\mathsf                              %  classes (of any things)
\let\oper=\mathbb                               %  operators
\bmdefine{\A}{A}                                %  particular algebras
\bmdefine{\2}{2}
\bmdefine{\B}{B}
\bmdefine{\D}{D}
\bmdefine{\M}{M}                                %  the monoid of substitutions
\bmdefine{\LLL}{L}                              %  algebraic language
\bmdefine{\Fm}{Fm}                              %  formula algebra
\bmdefine{\zerou}{[0{,}1]}  
\bmdefine{\T}{T}                                %  particular algebras
\newcommand{\VVV}{\oper{V}}                     %  class operators

\bmdefine{\boldstar}{\mathchoice{\textstyle*}{\textstyle*}{\textstyle*}{\scriptstyle*}}

\newcommand{\X}[0]{\mathbb{X}}

\bmdefine{\btau}{\tau}                                  %  transformer tau
\bmdefine{\brho}{\rho}                                  %  transformer rho
\newcommand{\down}[0]{{\downarrow}}
\newcommand{\up}[0]{{\uparrow}}

\bmdefine{\leibniz}{\mathbb{N}}        %  Leibniz operator

\bmdefine{\frege}{\Lambda}         %  Frege operator

\makeatletter
\newcommand{\tarskidsp}{\mathord%
   {\m@th\raisebox{0pt}[0pt][0pt]{$\stackrel%
   {\raisebox{-2.7pt}[0ex][0pt]{$\displaystyle \,\?\thicksim$}}%
   {\displaystyle\leibniz}$}}}
\newcommand{\tarskitxt}{\mathord%
   {\m@th\raisebox{0pt}[0pt][0pt]{$\stackrel%
   {\raisebox{-2.7pt}[0ex][0pt]{$\,\?\thicksim$}}{\displaystyle\leibniz}$}}}
\newcommand{\tarskiscr}{\mathord%
   {{\m@th\raisebox{0pt}[0pt][0pt]{$\stackrel%
   {\raisebox{-2.4pt}[0ex][0pt]{$\scriptstyle \,\?\thicksim$}}%
   {\scriptstyle\leibniz}$}}}}
\newcommand{\tarskiscrscr}{\mathord%
   {{\m@th\raisebox{0pt}[0pt][0pt]{$\stackrel%
   {\raisebox{-2pt}[0ex][0pt]{$\scriptscriptstyle \,\?\thicksim$}}%
   {\scriptscriptstyle\leibniz}$}}}}
\newcommand{\tarski}{\@ifnextchar ^ %
   {\mathchoice{\tarskidsp\kern-.07em}{\tarskitxt\kern-.07em}%
   {\tarskiscr\kern-.07em}{\tarskiscrscr\kern-.07em}}%
   {\mathchoice{\tarskidsp}{\tarskitxt}{\tarskiscr}{\tarskiscrscr}}}
\makeatother

\theoremstyle{theorem}
\newtheorem{Theorem}{Theorem}[section]
\newtheorem{Coloring Theorem}[Theorem]{Coloring Theorem}
\newtheorem{Claim}[Theorem]{Claim}
\newtheorem{Proposition}[Theorem]{Proposition}
\newtheorem{Lemma}[Theorem]{Lemma}
\newtheorem{Corollary}[Theorem]{Corollary}

\theoremstyle{definition}
\newtheorem{law}[Theorem]{Definition}

\theoremstyle{remark}

\newtheorem{Remark}[Theorem]{Remark}

\begin{document}
\title[On locally finite varieties of Heyting algebras]{On locally finite varieties of Heyting algebras}

\author{M. Martins and T. Moraschini}
\address{Departament de Filosofia, Facultat de Filosofia, Universitat de Barcelona (UB), Carrer Montalegre, $6$, $08001$ Barcelona, Spain}
%\email{jansana@ub.edu}
%\address{Institute of Computer Science, Academy of Sciences of Czech Republic, Pod Vod\'arenskou v\v{e}\v{z}\'{i} $271/2$, $182$ $07$ Prague $8$, Czech Republic}

\email{miguel.pedrosodelima@ub.edu \\
tommaso.moraschini@ub.edu}
\date{}

\maketitle

\begin{abstract}
For every $n \in \mathbb{N}$, we construct a variety of Heyting algebras, whose $n$-generated free algebra is finite but whose $(n+1)$-generated free algebra is infinite.
\end{abstract}

%Given a poset $\mathbb{X}$ and $x \in X$, sometimes we write ${\uparrow}^{\mathbb{X}}x$ to stress that this upset is computed in $\mathbb{X}$.

%Notation: given a map $f$, we denote by $\textup{Ker}(f)$ its kernel.

%Explain $\alpha$ and $\beta$-reductions.

\section{Introduction}

A \emph{Heyting algebra} is a structure $\langle A; \land, \lor, \to, 0, 1 \rangle$ where $\langle A; \land, \lor, 0, 1 \rangle$ is a bounded distributive lattice and $\to$ a binary operation satisfying \emph{residuation law}, that is, the demand that for every $a, b, c \in A$,
\[
a \land b \leq c \Longleftrightarrow a \leq b \to c.
\]
The class of Heyting algebras forms a \emph{variety} (i.e., an equational class) that we denote by $\mathsf{HA}$ \cite{BaDw74,ChZa97,Esakia-book85,RaSi63}.

From a logical standpoint, the importance of Heyting algebras is that they algebraize the \emph{intuitionistic propositional claculus} $\mathsf{IPC}$ in the sense of \cite{BP89}. As a consequence, the axiomatic extensions of $\mathsf{IPC}$ (known as \emph{superintuitionistic logics}, or si-logics for short) form a lattice that is dually isomorphic to that of varieties of Heyting algebras. This allows us to study each si-logic through the lenses of its corresponding variety of Heyting algebras which, in turn, is amenable to the methods of universal algebra and duality theory. 

As a particular instance of this phenomenon, we recall that an si-logic $\mathsf{L}$ is  \emph{locally tabular} when for every $n \in \mathbb{N}$ there are only finitely many formulas in variables $x_1, \dots, x_n$ up to logical equivalence in $\mathsf{L}$. On the other hand, a variety is called \emph{locally finite} when its finitely generated members (equiv.\ finitely generated free algebras) are finite \cite[Thm.\ II.10.15]{BuSa00}. It is well known that an si-logic is locally tabular iff the corresponding variety of Heyting algebras is locally finite.

A fascinating problem by Bezhanishvili and Grigolia asks to determine whether it is true that a variety of Heyting algebras is locally finite iff its two-generated free algebra is finite \cite[Prob.\ 2.4]{BezGri05}. While this holds in the restrictive context of varieties of Heyting algebras of width two \cite{TBenjam2020}, in this paper we establish the following:

\begin{Theorem}
For every $n\geq 2$ there exists a variety of Heyting algebras whose $n$-generated free algebra is finite, while its $(n+1)$-generated free algebra is infinite.
\end{Theorem}

\noindent Consequently, for every $n \geq 2$ there exits a nonlocally finite variety of Heyting algebra whose free $n$-generated algebra is finite.

This result was established in the spring of 2020, at a time when the second author was supervising the master thesis \cite{TBenjam2020}. Recently, an alternative proof was independently discovered by Hyttinen and Quadrellaro \cite{HyttQuad23}. This motivated us to share the original proof as well.

\section{Esakia spaces}

Let $\mathbb{X}= \langle X, \leq \rangle$ be a poset. We denote the \textit{upset generated} by a subset $U$ of $X$ by 
\[
\up U\coloneqq \{x\in X : \exists u\in U \text{ such that }u \leq x\},
\]
and if $U = \up U$, then $U$ is called an \textit{upset}. 
If $U=\{x\}$, we simply write ${\uparrow} x$ and call it a \textit{principal upset}. 
The notion of a \textit{downset} and the arrow operator $\down$ are defined analogously. 

If $x,y \in X$, then $x$ is said to be an \textit{immediate predecessor} of $y$ if $x < y$ and no point in $X$ lies between them (i.e., if $z\in X$ is such that $x \leq z \leq y$, then either $x=z$ or $y=z$). 
If this is the case, we call $y$ an \textit{immediate successor} of $x$. 

\begin{law} \label{def bi-esa}
A triple $\X= \langle X,\tau,\leq \rangle$ is an \textit{Esakia space} if it is a compact ordered topological space satisfying the following conditions:
\benroman 
    \item $\down U$ is clopen, for every clopen $U$;
    \item \textit{Priestley separation axiom}:  for every $x, y \in X$,
    \[
x \nleq y \text{ implies that there exists a clopen upset $U$ such that }x \in U \text{ and } y \notin U.
    \]
\eroman
\end{law}

\begin{law} \label{def p-morphism}
A  map $f\colon \mathbb{X} \to \mathbb{Y}$ between Esakia spaces is called an \textit{Esakia morphism} if it it continuous, order preserving, and for every $x \in X$ and $y \in Y$,
\[
\text{if }f(x) \leq y, \text{ there exists some } z\geq x \text{ such that }f(z)=y.
\]
\end{law}
When endowed with the discrete topology, every finite poset becomes an Esakia space. In fact, this is the only way to view a finite poset as an Esakia space, because Esakia spaces are Hausdorff.

In view of Esakia duality \cite{Es74,Esakia-book85}, the category $\mathsf{HA}$ of Heyting algebras with homomorphisms and the category $\mathsf{ES}$ of Esakia spaces with Esakia morphisms are dually equivalent. We denote the contravariant functors witnessing this duality by $(-)_\ast \colon \mathsf{HA} \to \mathsf{ES}$ and $(-)^\ast \colon \mathsf{ES} \to \mathsf{HA}$.

\begin{law}
Let $\mathbb{X}$ be an Esakia space.
\begin{enumerate}[\normalfont (1)]
\item An {\em E-subspace} of $\mathbb{X}$ is a closed upset equipped with the subspace topology and the restriction of the order.
\item An {\em E-partition} on $\mathbb{X}$ is an equivalence relation $R$ on $X$ such that for all $x, y, z \in X$:
\benroman
\item if $\langle x, y \rangle \in R$ and $x \leq z$, then $y\leq w$ and $\langle z, w \rangle \in R$ for some $w \in X$;
\item if $\langle x, y \rangle \notin R$, then there is an {\em $R$-saturated} clopen $U$ (a union of equivalence classes of $R$) such that $x \in U$ and $y \notin U$.
\eroman
\end{enumerate}
\end{law}

\begin{Remark}
    It follows from the definition of an Esakia space $\X$ that its principal upsets are closed.   Thus, $\up x$ can be viewed as an E-subspace of $\X$, for all $x \in X$.
\end{Remark}

If $R$ is an E-partition of the Esakia space $\X = \langle X, \tau, \leq \rangle$, we denote by $\mathbb{X} / R$ the Esakia space obtained by endowing $X / R$ with the quotient topology and the partial order $\sqsubseteq$ defined as follows:
\[
x / R \sqsubseteq y / R \Longleftrightarrow x' \leq y' \text{ for some }x' \in x / R \text{ and } y' \in y /R.
\]
Moreover, the natural map $f \colon \X \to \X/R$ is a surjective Esakia morphism. Conversely, the kernel $\mathsf{Ker}(f)$ of every surjective Esakia morphism $f \colon \X \to \mathbb{Y}$, is an E-partition of $\X$ such that $\X / \mathsf{Ker}(f) \cong \mathbb{Y}$ \cite[Cor.\ 2.3.1]{Bez-PhD}.

We will rely on the following observations \cite[Lems.\ 3.1.6 and 3.1.7]{Bez-PhD}:

\begin{Lemma}
    Let $\X$ be an Esakia space and $x,y \in X$.
    \benroman
        \item If $y$ is the unique immediate successor of $x$ and $R$ is the smallest equivalence relation that identifies $x$ and $y$, then $R$ is an E-partition. The natural map $f\colon \X \to \X/R$ is an Esakia morphism, called an $\alpha$\textit{-reduction.}

        \item If $x$ and $y$ have exactly the same immediate successors and $R$ is the smallest equivalence relation that identifies $x$ and $y$, then $R$ is an E-partition. The natural map $f\colon \X \to \X/R$ is an Esakia morphism, called a $\beta$\textit{-reduction.}
    \eroman
\end{Lemma}  

\begin{Lemma}
    If $f \colon \X \to \mathbb{Y}$ is an Esakia morphism between finite Esakia spaces, then there exists a finite sequence $f_1, \dots , f_n$ of $\alpha$ and $\beta$-reductions such that $f=f_n \circ \dots \circ f_1$.
\end{Lemma}

\section{Finitely generated Heyting algebras}

Finitely generated Heyting algebras can be described in terms of their Esakia duals, as we proceed to explain.  

\begin{law}
Given $n \in \mathbb{N}$, let $\mathbb{C}_{n} = \langle C_n; \sqsubseteq \rangle$ be the poset whose universe is the set of sequences of length $n$ of zeros and ones and whose order is defined as follows:
\[
\langle m_{1}, \dots, m_{n}\rangle \sqsubseteq \langle k_{1}, \dots, k_{n}\rangle \Longleftrightarrow m_i \leq k_i, \text{for every }i \leq n.
\]
\end{law}

The elements of $\mathbb{C}_n$ can be used to color Esakia spaces as follows.  Given an Esakia space $\mathbb{X}$ and a function $f \colon \mathbb{X} \to \mathbb{C}_n$, we say that an element $x \in X$ is \textit{colored} by a color $\vec{c} \in C_n$ when $f(x) = \vec{c}$. The set of elements of $\mathbb{X}$ colored by $\vec{c}$ will be denoted by $\vec{c}(\mathbb{X})$.

\begin{law}
Let $\mathbb{X}$ be an Esakia space and $n \in \mathbb{N}$. A function $f \colon \mathbb{X} \to \mathbb{C}_n$ is said to be
\benroman
\item a \textit{weak $n$-coloring} of $\mathbb{X}$ if $f$  is order preserving and $\vec{c}(\mathbb{X})$ is clopen, for all $\vec{c} \in C_{n}$;
\item an \textit{$n$-coloring} of $\mathbb{X}$ if it is a weak $n$-coloring such that every E-partition of $\mathbb{X}$ other than the identity relation identifies two elements of $\mathbb{X}$ of distinct color. 
\eroman
An Esakia space $\mathbb{X}$ is said to be $n$\textit{-colorable}  when there it admits an $n$-coloring.
\end{law}

Finitely generated Heyting algebras and colorings of Esakia spaces are related as follows \cite[Thm.\ 3.1.5]{Bez-PhD}. 

\begin{Coloring Theorem}\label{Thm:coloring} Let $\A$ be a Heyting algebra and $n \in \mathbb{N}$.  Then $\A$ is $n$-generated if and only if $\A_\ast$ is $n$-colorable.
\end{Coloring Theorem}

A Heyting algebra $\A$ is \emph{subdirectly irreducible} (SI for short) when it has a second greatest element. This is equivalent to the demand that $\A_*$ has a least element that, moreover, is isolated \cite[Prop.\ A.1.2]{Esakia-book85}. When a poset (or, in particular, an Esakia space) $\X$ has a least element $x$, we call $x$ the \textit{root} of $\X$, and say that $\X$ is \textit{rooted}.
Thus, an Heyting algebra $\A$ is SI iff $\A_*$ has an isolated root.

 In view of the following observation,  in order to determine whether a variety of Heyting algebras is locally finite, it suffices to examine its finite SI members \cite[Thm.\ 4.3]{BezMor19}.

\begin{Theorem}\label{Thm:lf vars}
A variety $\class{V}$ of Heyting algebras is locally finite if and only if $\class{V}$ has, up to isomorphism, only finitely many finite $n$-generated SI members, for every $n \in \mathbb{N}$.
\end{Theorem}

The following technical observation will be used throughout the paper.  In its statement, $\vec{0}$ stands for the color consisting of $n$ zeros.

\begin{Lemma}\label{Lem:local-finiteness}
Let $\{ \mathbb{X}_{m} : m \in \mathbb{N} \}$ be a family of finite $n$-colorable Esakia spaces such that
\[
\vert X_{1} \vert < \vert X_{2} \vert < \cdots < \vert X_{m} \vert < \cdots 
\]
and whose antichains are of size $\leq t$ for some $t \in \mathbb{N}$.\ Then there are $k \in \mathbb{N}$ and a family $\{ \mathbb{Z}_{m} : m \in \mathbb{N} \}$ of E-subspaces of spaces in $\{ \mathbb{X}_{m} : m \in \mathbb{N} \}$ such that 
\benroman
\item $\vert Z_{1} \vert < \vert Z_{2} \vert < \cdots < \vert Z_{m} \vert < \cdots$ and
\item each $\mathbb{Z}_{m}$ admits an $n$-coloring for which $\vert Z_{m} \smallsetminus \vec{0}(\mathbb{Z}_{m}) \vert \leq k$.
\eroman
\end{Lemma}

The proof of this lemma relies on the following easy combinatorial principle that follows, for instance, from Dilworth's Theorem \cite{Di50a}.

\begin{Proposition}\label{Prop:weak-Dilworth}
For every $n, m \in \mathbb{N}$, if chains and antichains in a poset $\mathbb{X}$ are, respectively, of size $\leq n$ and $\leq m$, then $\mathbb{X}$ has at most $n \times m$ elements.
\end{Proposition}

\begin{proof}[Proof of Lemma \ref{Lem:local-finiteness}]
Fix an $n$-coloring $f_m$ on each $\mathbb{X}_m$.  The consider the set of colors
\[
D \coloneqq \{ \vec{c} \in \{ 0, 1 \}^n : \text{for every }k \in \mathbb{N}\text{ there exists }m_k \in \mathbb{N} \text{ such that }k \leq \vert \vec{c}(\mathbb{X}_{m_k}) \vert \}.
\]
Notice that $D$ is nonempty,  because
\[
\vert X_{1} \vert < \vert X_{2} \vert < \cdots < \vert X_{m} \vert < \cdots
\]
Furthermore,  it is finite.  Therefore, when viewed as a subposet of $\mathbb{C}_n$,  the set $D$ has at least a maximal element $\vec{c}$.

For every $m \in \mathbb{N}$ and $x \in \vec{c}(\mathbb{X}_m)$,  let $\mathbb{Y}_{m}^{x}$ be the E-subspace of $\mathbb{X}_m$ with universe ${\uparrow}x$.  Moreover,  let $g_m^x \colon \mathbb{Y}_{m}^{x} \to \mathbb{C}_n$ be the function defined, for every $y \in {\uparrow}x$, as
\[
g_m^x(y) = 
\begin{cases}
f_m(y) & \text{ if } f_m(y) \ne \vec{c}\\
\vec{0} & \text{ otherwise.} 
\end{cases}
\]
We will use repeatedly the fact that $\vec{c} \sqsubseteq g_m^x(z)$, for every element $z \in Y_m^x$.  To prove this, observe that $x \leq z$, because $Y_m^x = {\uparrow}x$. Moreover,  $f_m(x) = \vec{c}$, since $x \in \vec{c}(\mathbb{X}_m)$. As $f_m$ is order preserving, we conclude that $\vec{c} = f_m(x) \sqsubseteq f_m(z)$, as desired.

\begin{Claim}\label{Claim:tech-1}
The map $g_m^x$ is an $n$-coloring of $\mathbb{Y}_{m}^{x}$.
\end{Claim}

\begin{proof}[Proof of the Claim.]
To prove that $g_m^x$ is order preserving, consider $y, z \in {\uparrow}x$ such that $y \leq z$.  If $f_m(y) = \vec{c}$, then from the definition of $g_m^x$ it follows
\[
g_m^x(y) = \vec{0} \sqsubseteq g_m^x(z),
\]
because $\vec{0}$ is the minimum of $\mathbb{C}_n$.\   Then we consider the case where $f_m(y) \ne \vec{c}$.  Since  $\vec{c} \sqsubseteq f_m(y)$, we get $\vec{c} \sqsubset f_m(y)$.  As $y \leq z$ and $f_m$ is order preserving, $\vec{c} \sqsubset f_{m}(z)$ also holds.  From the definition of $g_m^x$ it follows that $g_m^x(y) = f_m(y)$ and $f_m(z) = g_m^x(z)$. Since $y \leq z$ and $f_m$ is order preserving, we obtain that
\[
g_m^x(y) = f_m(y) \sqsubseteq f_m(z) = g_m^x(z).
\]
We conclude that $g_m^x$ is order preserving.\ Furthermore,  recall that $\mathbb{Y}_m^x$ is finite and, therefore,  its topology is discrete.\ As a consequence,  $\vec{a}(\mathbb{Y}_m^x)$ is clopen, for every color $\vec{a} \in \{ 0, 1 \}^n$.  Hence, $g_m^x$ is a weak $n$-coloring of $\mathbb{Y}_m^x$.

To prove that $g_m^x$ is also an $n$-coloring,  consider an E-partition $R$ of $\mathbb{Y}_m^x$ other than the identity relation on $Y_m^x$.  Moreover, let $id_{X_m}$ be the identity relation on $\mathbb{X}_m$.  We will prove that the union $S \coloneqq R \cup id_{X_m}$ is an E-partition on $\mathbb{X}_m$.  To this end, consider $y, z, w \in X_m$ such that $\langle y, z \rangle \in S$ and $y \leq w$.  We need to show that there exists $u \in X_m$ such that $u \geq z$ and $\langle w, u \rangle \in S$. If $y = z$, then we are done taking $u \coloneqq w$.  Then we consider the case where $y \ne z$. From the definition of $S$ it follows that $\langle y, z \rangle \in R$.  Since $y \leq w$ and $\mathbb{Y}_m^x$ is an upset of $\mathbb{X}_m$,  the element $w$ also belongs to $\mathbb{Y}_m^x$. Since $R$ is a an E-partition on $\mathbb{Y}_m^x$, there exists $u \geq z$ such that $\langle w, u \rangle \in R \subseteq S$. Since the topology of $\mathbb{X}_m$ is discrete (because $\mathbb{X}_m$ is finite), this shows that $S$ is an E-partition on $\mathbb{X}_m$.

Now,  recall that there exists a pair $\langle y, z \rangle \in R$ such that $y \ne z$.  As $R \subseteq S$,  the relation $S$ is an E-partition on $\mathbb{X}_m$ different from the identity relation $id_{X_m}$. Since $f_m$ is an $n$-coloring of $\mathbb{X}_m$,  there is a pair $\langle y^\ast, z^\ast\rangle \in S$ such that $f_m(y^\ast) \ne f_m(z^\ast)$.\  Clearly, $y^\ast \ne z^\ast$.\ As $S = R \cup id_{X_m}$, we obtain $\langle y^\ast, z^\ast \rangle \in R$.
 We will prove that $g_m^x(y^\ast) \ne g_m^x(z^\ast)$.  If both $f_m(y^\ast)$ and $f_m(z^\ast)$ are different from $\vec{c}$, the definition of $g_m^x$ implies that
 \[
 g_m^x(y^\ast) = f_m(y^\ast) \ne f_m(z^\ast) = g_m^x(z^\ast),
\]
as desired. Then we consider the case where $f_m(y^\ast)$ or $f_m(z^\ast)$ is $\vec{c}$. By symmetry, we may assume that $f_m(y^\ast) = \vec{c}$.  By the definition of $g_m^x$, this yields $g_m^x(y^\ast) = \vec{0}$.  Moreover,  as $f_m(y^\ast) \ne f_m(z^\ast)$, we have $f_m(z^\ast) \ne \vec{c}$ and, therefore, $\vec{c}\sqsubset f_m(z^\ast)$.  By the definition $g_m^x$, we obtain
\[
g_m^x(y^\ast) = \vec{0} \sqsubseteq \vec{c} \sqsubset f_m(z^\ast) = g_m^x(z^\ast).
\]
Therefore, $R$ identifies two distinct elements of $\mathbb{Y}_m^x$, namely $y^\ast$ and $z^\ast$,  colored differently by $g_m^x$.  Hence, we conclude that $g_m^x$ is an $n$-coloring of $\mathbb{Y}_m^x$.
\end{proof}

The sequence of posets $\mathbb{Y}_m^x$ with the colorings $g_m^x$ has the following property.

\begin{Claim}\label{Claim:tech-2}
There exists $k \in \mathbb{N}$ such that for all $m \in \mathbb{N}$ and $x \in \vec{c}(\mathbb{X}_m)$,
\[
\vert  \{ z \in Y_m^x : g_m^x(z) \ne \vec{0} \}\vert  \leq k.
\] 
\end{Claim}

\begin{proof}[Proof of the Claim.]
Let $M$ be the set of colors in $\{ 0, 1 \}^n$ strictly larger than $\vec{c}$ in $\mathbb{C}_n$.  Since $M$ is finite and, by assumption,  $\vec{c}$ is maximal in $D$, there exists $k \in \mathbb{N}$ such that
\[
\vert \{ z \in X_m : f_m(z) \in M \} \vert \leq k,
\]
for all $m \in \mathbb{N}$. Consequently, for every $m \in \mathbb{N}$ and $x \in \vec{c}(\mathbb{X}_m)$,
\begin{equation}\label{Eq:lemma-1}
\vert \{ z \in Y_m^x : g_m^x(z) \in M \} \vert \leq k.
\end{equation}

Now, consider $m\in \mathbb{N}$ and $x \in \vec{c}(\mathbb{X}_m)$.  We will prove that
\begin{equation}\label{Eq:lemma-2}
\{ z \in Y_m^x : g_m^x(z) \ne \vec{0} \} \subseteq \{ z \in Y_m^x : g_m^x(z) \in M \}.
\end{equation}
To this end, consider $z \in Y_m^x$ such that $g_m^x(z) \ne \vec{0}$.  By the definition of $g_m^x$, this implies $f_m(z) \ne \vec{c}$ and $g_m^x(z) = f_m(z)$. Since $\vec{c} \sqsubseteq f_m(z)$,  we obtain $\vec{c} \sqsubset f_m(z) = g_m^x(z)$. Hence,  we conclude that $g_m^x(z) \in M$. This establishes (\ref{Eq:lemma-2}). Lastly, from (\ref{Eq:lemma-1}) and (\ref{Eq:lemma-2}) it follows that $\vert  \{ z \in Y_m^x : g_m^x(z) \ne \vec{0} \} \vert \leq k$, as desired.
\end{proof}

\begin{Claim}\label{Claim:tech-3}
For every $k \in \mathbb{N}$ there are $m \in \mathbb{N}$ and $x \in \vec{c}(\mathbb{X}_{m})$ such that $k \leq \vert Y_{m}^x \vert$.
\end{Claim}

\begin{proof}[Proof of the Claim.]
Suppose the contrary with a view to contradiction.  Then there exists $k \in \mathbb{N}$ such that $\vert Y_{m}^x \vert \leq k$, for all $m \in \mathbb{N}$ and $x \in  \vec{c}(\mathbb{X}_{m})$.\ Thus, for each $m \in \mathbb{N}$,  the chains of the subposet $\vec{c}(\mathbb{X}_m)$ of $\mathbb{X}_m$ must have size $\leq k$.  Now, recall that antichains in $\mathbb{X}_m$ have size $\leq t$, by assumption. Therefore, antichains in $\vec{c}(\mathbb{X}_m)$ have also size $\leq t$. Using Proposition \ref{Prop:weak-Dilworth}, we conclude that $\vert \vec{c}(\mathbb{X}_m) \vert \leq k \times t$, for every $m \in \mathbb{N}$. But this contradicts the fact that $\vec{c}$ belongs to $D$.
\end{proof}

In view of Claims \ref{Claim:tech-2} and \ref{Claim:tech-3},  there exists a subset
\[
\{ \mathbb{Z}_m : m \in \mathbb{N} \} \subseteq \{ \mathbb{Y}_m^x : m \in \mathbb{N} \text{ and }x \in \vec{c}(\mathbb{X}_m) \}
\]
such that
\benroman
\item $\vert Z_{1} \vert < \vert Z_{2} \vert < \cdots < \vert Z_{m} \vert < \cdots$ (this is made possible by of Claim  \ref{Claim:tech-3} and the fact that the various $\mathbb{Y}_m^x$ are finite) and
\item there exists $k \in \mathbb{N}$ such that for every $m \in \mathbb{N}$,  
\[
\text{if }\mathbb{Z}_m = \mathbb{Y}_m^x\text{, then }\vert  \{ z \in Z_m : g_m^x(z) \ne \vec{0} \} \vert \leq k.
\]
(this is made possible by of Claim  \ref{Claim:tech-2}).
\eroman
As $g_m^x$ is an $n$-coloring of $\mathbb{Y}_m^x$ by Claim \ref{Claim:tech-1},  the sequence of posets $\{ \mathbb{Z}_m : m \in \mathbb{N} \}$ colored with the suitable $g_m^x$ (that is, if $\mathbb{Z}_m = \mathbb{Y}_m^x$, we color $\mathbb{Z}_m$ with $g_m^x$) satisfies the conditions in the statement.
\end{proof}

\section{Introducing the abominations}

Our first aim is to prove that for every $n \in \mathbb{N}$ there exists a variety $\class{K}_n$ of Heyting algebras  whose $n$-generated members are finite,  but whose $(n+1)$-generated ones may be infinite.  To this end, we will exhibit a special Esakia space $\mathbb{X}_n$, called the $n$-\textit{abomination},  and let $\class{K}_n$ be the variety generated by the algebraic dual of $\mathbb{X}_n$.

The next definition,  inspired by the one-point compactification of a discrete space, is instrumental in the construction of $\mathbb{X}_n$.

\begin{law}
 The \textit{root compactification} of a poset $\mathbb{X}$ is the ordered topological space obtained by adding a new minimum $\bot$ to $\mathbb{X}$ and declaring open a subset $U$ of $X \cup \{ \bot \}$ provided that 
\[
\text{if }\bot \in U \text{, then }U \text{ is cofinite}.
\]
\end{law}

\begin{Proposition}\label{Prop:root-compactification}
Let $\mathbb{X}$ be a poset such that ${\downarrow}x$ is cofinite for every $x \in X$. Then root compactification of $\mathbb{X}$ is an Esakia space.
\end{Proposition}

\begin{proof}
Let $\mathbb{X}_{\bot}$ be the root compactification of $\mathbb{X}$ with universe $X_{\bot} = X \cup \{ \bot \}$.  The fact that $\mathbb{X}_{\bot}$ is compact is an immediate consequence of its definition. To prove that it satisfies Priestley separation axiom, consider $x, y \in X_{\bot}$ such that $x \nleq y$.  By assumption the downset ${\downarrow}y$ of $y$ in $\mathbb{X}_\bot$ is cofinite.\ This guarantees that ${\downarrow}y$ is open. As $\bot \leq y$, the complement $X_\bot \smallsetminus {\downarrow}y$ does not contains $\bot$ and, therefore, is also open. It follows that ${\downarrow}y$ is a clopen downset of $\mathbb{X}_\bot$. that contains $y$ but omits $x$,  as desired.

To conclude that $\mathbb{X}_{\bot}$ is an Esakia space, it only remains to show that downsets of opens are open. To this end, consider an open set $U$.  If $U$ is empty, so is its downset and, since the empty set is open, we are done. Then we consider the case where $U$ contains some element $x$. By assumption, ${\downarrow}x$ is cofinite, whence  ${\downarrow}U$ is also cofinite.  This, in turn,  implies that ${\downarrow}U$ is open.
\end{proof}

%\begin{Theorem}
%For every $n \in \mathbb{N}$ there exists a nonlocally finite variety of Heyting algebras whose $n$-generated free algebra is finite.
%\end{Theorem}

Throughout this section, fix an integer $n \geq 2$. Then consider the set
\[
T_{n} \coloneqq \{ \langle k_{1}, k_{2}, k_{3}\rangle  : k_{1}, k_{2}, k_{3} \in \mathbb{N} \text{ are all different and}\leq 2^{n+1}-1\}
\]
and take an enumeration
\[
T_{n} = \{ s_{0}, \dots, s_{t}\}.
\]
Then consider disjoint sets of distinct elements 
\begin{align*}
A &\coloneqq \{ a_{m}  :  m \in \mathbb{N} \} \\
B  &\coloneqq \{ b_{m}  :  m \in \mathbb{N} \} \\
C  &\coloneqq\{ c_{m, k}  : k, m \in \mathbb{N} \text{ and } k \leq 2^{n+1}-1\} \\
D &\coloneqq \{ d_{m, k}  : k, m \in \mathbb{N} \text{ and } k \leq 2^{n+1}-1\} \\
 E^a &\coloneqq\{ e_{m, k}^{a}  : k, m \in \mathbb{N} \text{ and } k \leq 2^{n+1}-1\} \\
E^b &\coloneqq\{ e_{m, k}^{b}  : k, m \in \mathbb{N} \text{ and } k \leq 2^{n+1}-1\}
\end{align*}
and let
\[
U_{n} \coloneqq A \cup B \cup C \cup D \cup E^a \cup E^b.
\]
Moreover,  let $\prec$ be the binary relation on $U_{n}$ that relates two elements $x, y \in U_{n}$, in symbols $x \prec y$, if and only if one of the following conditions holds:
\benroman
\item\label{Eq:notstrange1} There is some $s_{j} = \langle k_{1}, k_{2}, k_{3}\rangle$ in $T_n$ and $m \in \mathbb{N}$ such that $m \equiv j \!\!\mod t+1$ and
\[
x = a_{m} \qquad y \in \{ c_{m,k_{1}}, c_{m,k_{2}} \} \qquad s_{j} = \langle k_{1}, k_{2}, k_{3}\rangle;
\]
\item\label{Eq:notstrange2} There is some $s_{j} = \langle k_{1}, k_{2}, k_{3}\rangle$ in $T_n$ and $m \in \mathbb{N}$ such that $m \equiv j \!\!\mod t+1$ and
\[
x = b_{m} \qquad y \in \{ c_{m,k_{1}}, c_{m,k_{3}} \} \qquad s_{j} = \langle k_{1}, k_{2}, k_{3}\rangle;
\]
\item\label{Eq:strange} There are $m, k, j, i\in \mathbb{N}$ such that $1 \leq m$ and $k \ne j$ and 
\[
x = c_{m, k} \qquad y \in \{ e_{m-1, j}^{a}, e_{m-1, i}^{b}\};
\]
\item\label{Eq:notstrange4}  There are $m, k, j \in \mathbb{N}$ such that $k \ne j$ and
\[
x = d_{m, k} \qquad y = c_{m, j};
\]
\item\label{Eq:notstrange5} There are $m, k, j \in \mathbb{N}$ such that $k \ne j$ and 
\[
x = e_{m, k}^{a} \qquad y \in \{ d_{m, j}, a_{m}\};
\]
\item\label{Eq:notstrange6} There are $m, k, j \in \mathbb{N}$ such that $k \ne j$ and 
\[
x = e_{m, k}^{b} \qquad y \in \{ d_{m, j}, b_{m}\}.
\]
\eroman
\noindent Notice that in Condition (\ref{Eq:strange}) is asymmetric, in the sense that the integer $j$ is assumed to be different from $k$, while $i$ might be equal $k$.

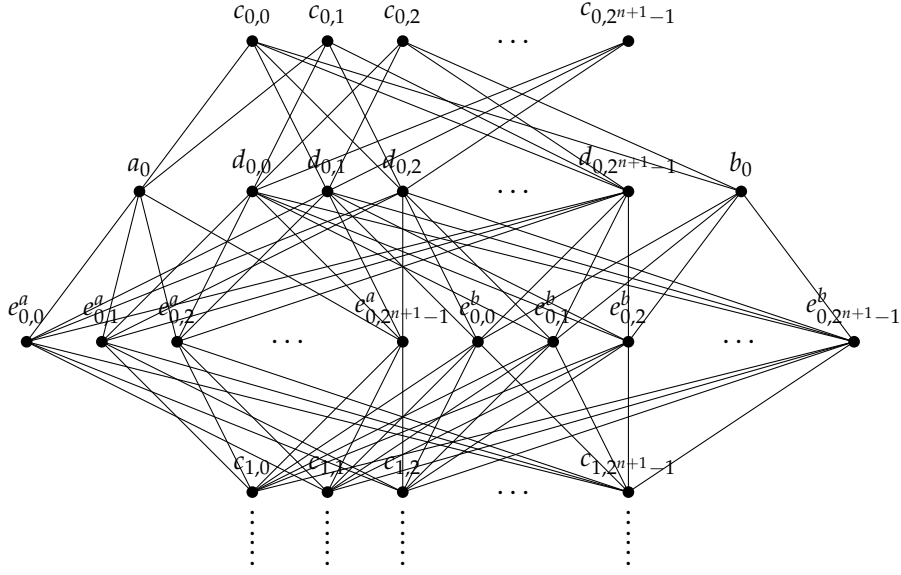
\begin{figure}[h]
\begin{tabular}{ccccccc}
  %\resizebox{0.5\textwidth}{!}{\begin{minipage}{\textwidth}
\\
\begin{tikzpicture}[node distance=9cm and 1cm,auto]
    \tikzstyle{point} = [shape=circle, thick, draw=black, fill=black , scale=0.35]
    \node  (c0) at (-2.5,2) [point,label={\small$c_{0,0}$}] {};
        \node  (c1) at (-1.5,2) [point,label={\small$c_{0,1}$}] {};
        \node  (c2) at (-0.5,2) [point,label={\small$c_{0,2}$}] {};
        \node  (cn) at (2.5,2) [point,label={\small$c_{0,2^{n+1}-1}$}] {};

 \node (zz1) at (1,2) {$\dots$};
  \node (zz2) at (1,0) {$\dots$};
 \node (zz3) at (-2,-2) {$\dots$};
  \node (zz3b) at (4,-2) {$\dots$};

 \node (zz4) at (1,-4) {$\dots$};

     \node  (d0) at (-2.5,0) [point,label={\small$d_{0,0}$}] {};
     \node  (d1) at (-1.5,0) [point,label={\small$d_{0,1}$}] {};
     \node  (d2) at (-0.5,0) [point,label={\small$d_{0,2}$}] {};
     \node  (dn) at (2.5,0) [point,label={\small$d_{0,2^{n+1}-1}$}] {};
     
     \node (a0) at (-4,0) [point,label={\small$a_{0}$}]  {};
     \node (b0) at (4,0) [point,label={\small$b_{0}$}]  {};

     \node  (ea0) at (-5.5,-2) [point,label={\small$e_{0,0}^{a}$}] {};
     \node  (ea1) at (-4.5,-2) [point,label={\small$e_{0,1}^{a}$}] {};
     \node  (ea2) at (-3.5,-2) [point,label={\small$e_{0,2}^{a}$}] {};
     \node  (ean) at (-0.5,-2) [point,label={\small$e_{0,2^{n+1}-1}^{a}$}] {};
     
     \node  (eb0) at (0.5,-2) [point,label={\small$e_{0,0}^{b}$}] {};
     \node  (eb1) at (1.5,-2) [point,label={\small$e_{0,1}^{b}$}] {};
     \node  (eb2) at (2.5,-2) [point,label={\small$e_{0,2}^{b}$}] {};
     \node  (ebn) at (5.5,-2) [point,label={\small$e_{0,2^{n+1}-1}^{b}$}] {};

     \node  (c10) at (-2.5,-4) [point,label={\small$c_{1,0}$}] {};
        \node  (c11) at (-1.5,-4) [point,label={\small$c_{1,1}$}] {};
        \node  (c12) at (-0.5,-4) [point,label={\small$c_{1,2}$}] {};
        \node  (c1n) at (2.5,-4) [point,label={\small$c_{1,2^{n+1}-1}$}] {};
     
    \node (dots0a) at (-2.5,-4.3) {$\vdots$};
    \node (dots0b) at (-2.5,-4.7) {$\vdots$};
     \node (dots1a) at (-1.5,-4.3) {$\vdots$};
    \node (dots1b) at (-1.5,-4.7) {$\vdots$};
    \node (dots2a) at (-0.5,-4.3) {$\vdots$};
    \node (dots2b) at (-0.5,-4.7) {$\vdots$};
    \node (dots3a) at (2.5,-4.3) {$\vdots$};    
    \node (dots3b) at (2.5,-4.7) {$\vdots$};
%    
%    

 %     \node (dots1) at (1,-5) {$\vdots$};
 %  \node (dots2) at (1,-5.4) {$\vdots$};
%\node (dots4) at (1,-6.3)  [point,label={\small$\bot$}] {};

 \draw  (a0) -- (c0)  ;
  \draw  (a0) -- (c1)  ;
  \draw  (b0) -- (c0)  ;
  \draw  (b0) -- (c2)  ;

    \draw  (d0) -- (c1)  ;
  \draw   (d0) -- (c2);
  \draw   (d0) -- (cn) ;
  \draw  (d1) -- (c0)  ;
  \draw   (d1) -- (c2);
  \draw   (d1) -- (cn) ;
  \draw  (d2) -- (c0)  ;
  \draw   (d2) -- (c1);
  \draw   (d2) -- (cn) ;
  \draw  (dn) -- (c0)  ;
  \draw   (dn) -- (c1);
  \draw   (dn) -- (c2) ;
  
    \draw   (ea0) -- (a0) ;
      \draw   (ea1) -- (a0) ;
          \draw   (ea2) -- (a0) ;
              \draw   (ean) -- (a0) ;
              
  \draw  (ea0) -- (d1)  ;
  \draw   (ea0) -- (d2);
  \draw   (ea0) -- (dn) ;
  \draw  (ea1) -- (d0)  ;
  \draw   (ea1) -- (d2);
  \draw   (ea1) -- (dn) ;
  \draw  (ea2) -- (d0)  ;
  \draw   (ea2) -- (d1);
  \draw   (ea2) -- (dn) ;
  \draw  (ean) -- (d0)  ;
  \draw   (ean) -- (d1);
  \draw   (ean) -- (d2) ;
  
  \draw   (eb0) -- (b0) ;
      \draw   (eb1) -- (b0) ;
          \draw   (eb2) -- (b0) ;
              \draw   (ebn) -- (b0) ;
              
  \draw  (eb0) -- (d1)  ;
  \draw   (eb0) -- (d2);
  \draw   (eb0) -- (dn) ;
  \draw  (eb1) -- (d0)  ;
  \draw   (eb1) -- (d2);
  \draw   (eb1) -- (dn) ;
  \draw  (eb2) -- (d0)  ;
  \draw   (eb2) -- (d1);
  \draw   (eb2) -- (dn) ;
  \draw  (ebn) -- (d0)  ;
  \draw   (ebn) -- (d1);
  \draw   (ebn) -- (d2) ;
  
   \draw  (c10) -- (ea1)  ;
  \draw   (c10) -- (ea2);
  \draw   (c10) -- (ean) ;
  \draw  (c11) -- (ea0)  ;
  \draw   (c11) -- (ea2);
  \draw   (c11) -- (ean) ;
  \draw  (c12) -- (ea0)  ;
  \draw   (c12) -- (ea1);
  \draw   (c12) -- (ean) ;
  \draw  (c1n) -- (ea0)  ;
  \draw   (c1n) -- (ea1);
  \draw   (c1n) -- (ea2) ;
  
      \draw  (c10) -- (eb0)  ;
    \draw  (c10) -- (eb1)  ;
  \draw   (c10) -- (eb2);
  \draw   (c10) -- (ebn) ;
  \draw  (c11) -- (eb0)  ;
      \draw  (c11) -- (eb1)  ;
  \draw   (c11) -- (eb2);
  \draw   (c11) -- (ebn) ;
  \draw  (c12) -- (eb0)  ;
     \draw  (c12) -- (eb2)  ;
  \draw   (c12) -- (eb1);
  \draw   (c12) -- (ebn) ;
  \draw  (c1n) -- (eb0)  ;
  \draw   (c1n) -- (eb1);
  \draw   (c1n) -- (eb2) ;
      \draw  (c1n) -- (ebn)  ;

      % \draw  (11) -- (00)  (11) -- (02) (11) -- (0n) 

\end{tikzpicture}
\end{tabular}
\caption{The top part of the reflexive and transitive closure of $\prec$ on $U_{n}$.}
\label{Fig:Xn}
\end{figure}

As an exemplification, the top part of the reflexive and transitive closure of $\prec$ on $U_{n}$ is depicted in Figure \ref{Fig:Xn}, under the assumption that $s_{0} = \langle 0, 1, 2 \rangle$.

\begin{law}
The \textit{abomination} $\mathbb{X}_{n}$ is the root compactification of the poset obtained by endowing $U_{n}$ with the reflexive and transitive closure of $\prec$.
\end{law}
%Let $\mathbb{X}_{n}$ be the poset obtained adding a minimum element $\bot$ to the direct union
%\[
%\bigcup_{m \in \mathbb{N}} \mathbb{X}_{m}^{n}.
%\]
\noindent The \textit{covering relation} of a poset $\mathbb{Z}$ is 
\[
R \coloneqq \{ \langle x, y \rangle \in Z \times Z : x\text{ is an immediate predecessor of }y \}.
\]
Notice that the covering relation of $\mathbb{X}_{n}$ is precisely $\prec$.

\begin{Proposition} \label{abomination is esakia}
For every integer $n \geq 2$,  the abomination $\mathbb{X}_{n}$ is an Esakia space. 
\end{Proposition}

\begin{proof}
Let $\leq$ be the reflexive and transitive closure of $\prec$ on $U_n$.\ We will rely on the following property of $\mathbb{X}_n$.

\begin{Claim}\label{Claim:X-is-Esakia1}
For all $m, k, j\in \mathbb{N}$ such that $k, j \leq 2^{n+1}-1$,
\[
\{ a_{m+1}, b_{m+1}, c_{m+1, j}, d_{m+1, j} e_{m+1, j}^a,  e_{m+1, j}^b\} \subseteq {\downarrow} e_{m, k}^b.
\]
\end{Claim}

\begin{proof}[Proof of the Claim.]
Consider $m, k \in \mathbb{N}$ such that $k \leq 2^{n+1}-1$.  By Condition (\ref{Eq:strange}) in the definition of $\prec$,  we obtain
\begin{equation}\label{Eq:X-is-Esakia-space1}
c_{m+1, j} \prec e_{m, k}^b \text{, for all }j \leq 2^{n+1}-1.
\end{equation}
Moreover,  from Conditions (\ref{Eq:notstrange1}), (\ref{Eq:notstrange2}), (\ref{Eq:notstrange4}) in the definition of $\prec$ it follows that
\begin{align}\label{Eq:X-is-Esakia-space2}
    \begin{split}
  & \text{for all }j \leq 2^{n+1}-1\text{ there are }p, q, r  \leq 2^{n+1}-1\\
   \text{ s.t.  }&a_{m+1} \prec c_{m+1, p}\text{ and }b_{m+1} \prec c_{m+1, q} \text{ and }d_{m+1, j} \prec c_{m+1, r}.
\end{split}
\end{align}
Lastly,  by Conditions (\ref{Eq:notstrange5}) and (\ref{Eq:notstrange6}) in the definition of $\prec$, we obtain
\begin{equation}\label{Eq:X-is-Esakia-space3}
e_{m+1, j}^a \prec a_{m+1}\text{ and }e_{m+1, j}^b \prec b_{m+1} \text{, for all }j \leq 2^{n+1}-1.
\end{equation}
Since $\leq$ is the reflexive and transitive closure of $\prec$, the result follows immediately from Conditions (\ref{Eq:X-is-Esakia-space1}),  (\ref{Eq:X-is-Esakia-space2}) and (\ref{Eq:X-is-Esakia-space3}).
\end{proof}

We will use the Claim to show that principal downsets in $\langle U_m; \leq \rangle$ are cofinite.  Indeed, by the Claim, this is true for principal downsets of the form ${\downarrow} e_{m, k}^b$. Therefore, it suffices to show that for all $x \in U_n$ there exist $m, k \in \mathbb{N}$ such that $e_{m, k}^b \leq x$.  If the element $x$ has already the form $e_{m, k}^b$, this is obvious.  Otherwise,  one of the the following cases holds for some $m, k \in \mathbb{N}$:
\benormal
\item\label{x:case2} $x$ has the form $b_m$;
\item\label{x:case4} $x$ has the form $d_{m, k}$;
\item\label{x:case3} $x$ has the form $c_{m, k}$;
\item\label{x:case5} $x$ has the form $e_{m,k}^a$;
\item\label{x:case1} $x$ has the form $a_m$.

\enormal

(\ref{x:case2}): Condition (\ref{Eq:notstrange6}) in the definition of $\prec$ implies that $e_{j,k}^b \leq b_m = x$ for all $j \leq 2^{n+1}-1$.

(\ref{x:case4}): Let $j \leq 2^{n+1}-1$ be distinct from $k$.  Condition (\ref{Eq:notstrange6}) in the definition of $\prec$ implies that $e_{m,j}^b \leq d_{m,k} = x$.

(\ref{x:case3}): Similarly, let $j \leq 2^{n+1}-1$ be distinct from $k$.  Condition (\ref{Eq:notstrange4}) in the definition of $\prec$ implies that $d_{m,j} \leq c_{m, k}$.  Therefore,  the result follows from case (\ref{x:case4}).

(\ref{x:case5}): Condition (\ref{Eq:strange}) in the definition of $\prec$ entails that there exists $j \leq 2^{n+1} - 1$ distinct from $k$ such that $c_{m+1,j} \leq e_{m,k}^a$. Therefore, the result follows from case (\ref{x:case3}).

(\ref{x:case1}): Condition (\ref{Eq:notstrange5}) in the definition of $\prec$ implies that $e_{m,k}^a \leq a_m $. Therefore,  the result follows from case (\ref{x:case5}).

\noindent Hence, we conclude that principal downsets in $\langle U_m; \leq \rangle$ are cofinite, as desired.  Together with Proposition \ref{Prop:root-compactification}, this implies that $\mathbb{X}_n$ is an Esakia space.  
\end{proof}

In order to keep proofs at a reasonable size,  from now on we will assume a basic understanding of the order structure of the abominations and, consequently,  omit precise references to the conditions in the definition of $\prec$ in our arguments.\ We believe, however, that the reader might find helpful to consult Figure \ref{Fig:Xn} when reading them.

Let $n$ be a positive integer. An element $x$ of a poset $\mathbb{X}$ has \emph{depth} $\leq n$ if every chain in $\up x$ has at most $n$ elements. If, moreover, $\up x$ contains a chain of $n$ elements, we say that $x$ has \emph{depth} $n$. Similarly, a rooted poset $\mathbb{X}$ is said to have \textit{width} $\leq n$ when all the antichains in $\mathbb{X}$ have at most $n$ elements. If, moreover, $\mathbb{X}$ has an antichain with $n$ elements, we say that $\mathbb{X}$ has \emph{width} $n$. An arbitrary (possibly nonrooted) poset has \emph{width }$\leq n$ (resp.\  $n$) when so do all its principal upsets. The empty poset is assumed to have width $0$.

%By construction, the width of $n$-abominations is bounded, as we proceed to explain.

\begin{Proposition} \label{width of abomination}
For every integer $n \geq 2$,  the $n$-abomination $\mathbb{X}_{n}$ is of width $2^{n+2}$.
\end{Proposition}

\begin{proof}

First, observe that for every $m \in \mathbb{N}$ the following is an antichain of size $2^{n+2}$ in $\mathbb{X}_n$:
\[
\{ e_{m,k}^{a} : 2^{n+1}-1 \geq k \in \mathbb{N} \} \cup \{ e_{m,k}^{b} : 2^{n+1}-1 \geq k \in \mathbb{N} \}.
\]
Since $\mathbb{X}_n$ is rooted, this implies that $\mathbb{X}_n$ cannot have width $\leq m$ for some $m < 2^{n+2}$.  Therefore, to conclude the proof, it suffices to show that antichains in $\mathbb{X}_n$ have size at most $2^{n+2}$. To this end, we rely on the following property of antichains in $\mathbb{X}_n$.

\begin{Claim}
For every antichain $Y$ in $\mathbb{X}_n$,
\begin{equation}
\text{if }Y \cap \{ a_m, b_m, c_{m, k}, d_{m, k} \} \ne \emptyset \text{ for some }m, k \in \mathbb{N}, \text{ then }\vert Y \vert \leq 2^{n+2}.
\end{equation}
\end{Claim}

\begin{proof}[Proof of the Claim.]
We have four cases, depending of whether $a_m, b_m, c_{m, k}$, or $d_{m, k}$ belongs to $Y$. Suppose first that $d_{m, k} \in Y$. The definition of $\prec$ implies that
\begin{align}\label{Eq:Y-antichain-dm}
X_n \smallsetminus ({\uparrow}d_{m, k} \cup {\downarrow}d_{m, k}) =& \, \,  \{ a_{m}, b_{m}, c_{m, k}, e^{a}_{m, k}, e^{b}_{m, k} \} \cup \{ d_{m, j} : k \ne j \leq 2^{n+1}-1\}.
\end{align}
Moreover, since $d_{m, k} \in Y$ and $Y$ is an antichain,
\[
Y \subseteq \{ d_{m,k} \} \cup (X_n \smallsetminus ({\uparrow}d_{m,k} \cup {\downarrow}d_{m,k})).
\]
Together with Conditions (\ref{Eq:Y-antichain-dm}), this implies $\vert Y \vert \leq 2^{n+1} + 5$. As, by assumption $n \geq 2$, we obtain $5 \leq 2^{n+1}$ and, therefore, 
\[
\vert Y \vert \leq 2^{n+1} + 5 \leq 2^{n+1} + 2^{n+1} = 2^{n+2}.
\]

Consequently, we may assume, without loss of generality, that $Y$ is disjoint from the set $\{ d_{m,j} : j \leq  2^{n+1}-1 \}$.  Then we consider the case where $c_{m,k} \in Y$.  The definition of $\prec$ implies that
\[
X_n \smallsetminus ({\uparrow}c_{m, k} \cup {\downarrow}c_{m, k}) \subseteq \, \,  \{ a_{m}, b_m, d_{m, k}, e^{a}_{m-1, k} \} \cup \{ c_{m, j} : k \ne j \leq 2^{n+1}-1\}.
\]
Moreover, since $c_{m, k} \in Y$ and $Y$ is an antichain, $Y \subseteq \{ c_{m,k} \} \cup (X_n \smallsetminus ({\uparrow}c_{m,k} \cup {\downarrow}c_{m,k}))$. Together with the above display, this implies that $\vert Y \vert \leq 2^{n+1}+4 \leq 2^{n+2}$, as desired.

As a consequence, we may assume, without loss of generality, that $Y$ is disjoint from the set $\{ c_{m,j} : j \leq  2^{n+1}-1 \}$.  Then we consider the case where $b_m \in Y$.  The definition of $\prec$ implies that
\begin{align}\label{Eq:Y-antichain-bm}
 \begin{split}
X_n \smallsetminus ({\uparrow}b_m \cup {\downarrow}b_m) =& \, \,  \{ e^{a}_{m, k} : k \leq 2^{n+1}-1\} \cup \{ d_{m, k} : k \leq 2^{n+1}-1\} \cup \{ a_m \} \cup\\
& \, \, \{ c_{m, k} : k \leq 2^{n+1}-1 \text{ and }k \notin \{ k_1, k_3 \} \},
\end{split}
\end{align}
where $k_1$ and $k_3$ are the unique positive integers with the following property: the unique positive integer $j \leq 2^{n+1}-1$ such that $m \equiv j \!\!\mod t+1$ is such that $s_j = \langle k_1, k_2, k_3 \rangle$, for some positive integer $k_2$.  Moreover, since $b_m \in Y$ and $Y$ is an antichain,
\[
Y \subseteq \{ b_m \} \cup (X_n \smallsetminus ({\uparrow}b_m \cup {\downarrow}b_m)).
\]
Together with Condition (\ref{Eq:Y-antichain-bm}) and the fact that $Y$ does not contain elements of the form $c_{m, k}$ or $d_{m, k}$, this yields 
\[
Y \subseteq \{ a_{m}, b_{m} \} \cup  \{ e^{a}_{m, k} : k \leq 2^{n+1}-1\}.
\]
Hence, we conclude that $\vert Y \vert \leq 2^{n+1}+2 \leq 2^{n+2}$. 

Lastly, we consider the case where $a_m \in Y$. The definition of $\prec$ implies that
\begin{align}\label{Eq:Y-antichain-am}
 \begin{split}
X_n \smallsetminus ({\uparrow}a_m \cup {\downarrow}a_m) =& \, \,  \{ e^{b}_{m, k} : k \leq 2^{n+1}-1\} \cup \{ d_{m, k} : k \leq 2^{n+1}-1\} \cup \{ b_m \} \cup\\
& \, \, \{ c_{m, k} : k \leq 2^{n+1}-1 \text{ and }k \notin \{ k_1, k_2 \} \},
\end{split}
\end{align}
where $k_1$ and $k_2$ are the unique positive integers with the following property: the unique positive integer $j \leq 2^{n+1}-1$ such that $m \equiv j \!\!\mod t+1$ is such that $s_j = \langle k_1, k_2, k_3 \rangle$, for some positive integer $k_3$. Moreover, since $a_m \in Y$ and $Y$ is an antichain,
\[
Y \subseteq \{ a_m \} \cup (X_n \smallsetminus ({\uparrow}a_m \cup {\downarrow}a_m)).
\]
Together with Condition (\ref{Eq:Y-antichain-bm}) and the fact that $Y$ does not contain elements of the form $c_{m, k}$ or $d_{m, k}$, this yields 
\[
Y \subseteq \{ a_{m}, b_{m} \} \cup  \{ e^{b}_{m, k} : k \leq 2^{n+1}-1\}.
\]
Hence, $\vert Y \vert \leq 2^{n+1}+2 \leq 2^{n+2}$. This concludes the proof of the Claim.
\end{proof}

We are now ready to prove that every antichain $Y$ in $\mathbb{X}_n$ has cardinality at most $2^{n+2}$.  First, the Claim guarantees that if $Y$ contains an element of the form $a_m, b_m, c_{m, k}$,  or $d_{m, k}$, then $\vert Y \vert \leq 2^{n+2}$.  Furthermore,  if $Y$ contains the minimum $\bot$ of $\mathbb{X}_n$,  then necessarily $Y = \{ \bot \}$ and, therefore, $\vert Y \vert = 1 \leq 2^{n+2}$.  Because of this, we may assume that
\[
Y \subseteq \{ e^{a}_{m,k} : m \in \mathbb{N} \text{ and } k \leq 2^{n+1}-1 \} \cup \{ e^{b}_{m,k} : m \in \mathbb{N} \text{ and } k \leq 2^{n+1}-1 \}.
\]
Let us suppose, with a view to contradiction, that $e_{m,i}^a,e_{m',j}^a\in Y$, for some $m \neq m' \in \mathbb{N}$ and $i,j \leq 2^{n+1}-1$.
Without loss of generality, we assume $m < m'$.
As $n \geq 2$, we know $2^{n+1}-1 \geq 7$, so there are $k \neq k' \leq 2^{n+1}-1$, both distinct from $i $ and $j$.
It follows from the definition of $\X_n$ that 
\[
e_{m',j}^a \leq d_{m',k} \leq c_{m',k'} \leq e_{m'-1,i}^a \leq d_{m'-1,k} \leq c_{m'-1,k'}\leq \dots \leq e_{m,i}^a,
\]
contradicting the assumption that $Y$ is an antichain.
An analogous argument shows that if $e_{m,i}^b,e_{m',j}^b\in Y$ for some $m, m' \in \mathbb{N}$ and $i,j \leq 2^{n+1}-1$, then we must have $m=m'$.

Suppose now that $e_{m,i}^a,e_{m',j}^b\in Y$, for some $m \neq m' \in \mathbb{N}$ and $i,j \leq 2^{n+1}-1$. 
If $m < m'$, then the definition of $\X_n$ entails
\[
e_{m',j}^b \leq d_{m',k} \leq c_{m',k'} \leq e_{m'-1,i}^a,
\]
for some $k \neq k' \leq 2^{n+1}-1$, both distinct from $i $ and $j$.
If $m=m'-1$, then the above display already contradicts the assumption that $Y$ is an antichain. If not, then we repeat our previous argument, again contradicting the aforementioned assumption.

If $m' <m $, we instead use the following inequality to derive a contradiction in a similar manner as above
\[
e_{m,i}^a \leq d_{m,k} \leq c_{m,k'} \leq e_{m-1,j}^b.
\]

From this discussion, we conclude that if an antichain $Y$ is contained in 
\[
 \{ e^{a}_{m,k} : m \in \mathbb{N} \text{ and } k \leq 2^{n+1}-1 \} \cup \{ e^{b}_{m,k} : m \in \mathbb{N} \text{ and } k \leq 2^{n+1}-1 \}.
\]
then $Y$ must be contained in 
\[
\{ e_{m,k}^{a} : 2^{n+1}-1 \geq k \in \mathbb{N} \} \cup \{ e_{m,k}^{b} : 2^{n+1}-1 \geq k \in \mathbb{N} \}.
\]
for some $m \in \mathbb{N}$, and therefore $|Y| \leq 2^{n+2}$, as desired.
\end{proof}

\begin{Proposition}\label{Prop:n+1-colorable}
$\X_{n}$ is $(n+1)$-colorable, for every integer $n \geq 2$.
\end{Proposition}

\begin{proof}
%We begin by showing that $\X_{n}$ is an Esakia space. Clearly, $\X$ is compact. On the other hand, it is easy to see that $\X$ satisfies Priestley's separation axiom. As downsets of opens are obviously open in $\X$, we conclude that $\X$ is an Esakia space. To prove that it satisfies Priestley separation axiom, consider $x, y \in X_{n}$ such that $x \nleq y$. Then $x \ne \bot$. Consequently, ${\uparrow}x$ is an upset containing $x$ and omitting $\bot$ and $y$. Since ${\uparrow}x$ omits $\bot$, by definition of the topology, $U$ is open. Furthermore, by Proposition \ref{Lem:construction}(\ref{item:construction1}), ${\uparrow}x$ is finite, whence $X_{n} \smallsetminus {\uparrow}x$ is cofinite. By definition of the topology, $X_{n} \smallsetminus {\uparrow}x$ is open, whence ${\uparrow}$ is a clopen upset that contains $x$ and omits $y$. Thus, $\X_{n}$ is a Priestley space. To conclude that $\X_{n}$ is an Esakia space, it only remains to show that downsets of opens are open. To this end, consider an open set $U$. We can assume without loss of generality that $U$ is nonvoid. By Proposition \ref{Lem:construction}(\ref{item:construction1}), ${\downarrow}U$ is cofinite and, therefore, open. We conclude that $\X_{n}$ is an Esakia space.
%
To prove that $\X_{n}$ is $(n+1)$-colorable, notice that it has $2^{n+1}$ maximal elements, namely
\[
c_{0, 0}, \dots, c_{0, 2^{n+1}-1}.
\]
Having at our disposal $2^{n+1}$ colors, we can color each maximal with a distinct color. Lastly, we color all the nonmaximal elements by the constant sequence $\vec{0}$. Because of the definition of the topology and order of $\X_{n}$, this is a weak $(n+1)$-coloring of $\X_{n}$.

To prove that it is an $(n+1)$-coloring, it only remains to show that every E-partition $R$ on $\X_{n}$ other than the identity relation identifies a pair of elements with distinct color. Suppose the contrary, with a view to contradiction, i.e., that there is an E-partition $R$ on $\X_{n}$ distinct from the identity that does not identify any pair of elements of the distinct color. Observe that there are distinct $\hat{x}, \hat{y} \in X_{n} \smallsetminus \{ \bot \}$ such that $\langle \hat{x}, \hat{y} \rangle \in R$. To prove this, recall that there are distinct $x$ and $y$ such that $\langle x, y \rangle \in R$.  If $x$ and $y$ are different from $\bot$, we are done. Then suppose that $y = \bot$, in which case $x \ne \bot$. By the definition of the order of $\X_{n}$, there is an element $\hat{y} > \bot$ such that $x \nleq \hat{y}$. Since $R$ is an E-partition and $\langle x, \bot \rangle \in R$, there is $\hat{x} \geq x$ such that $\langle \hat{x}, \hat{y}\rangle \in R$. By construction, $\hat{x}$ and $\hat{y}$ are different and other than $\bot$, as desired.

Recall that the proof of Proposition \ref{abomination is esakia} established $\X_n$ as an Esakia space whose principal downsets are cofinite. This, together with the fact that $\hat{x}$ and $\hat{y}$ are different from $\bot$, entails that the upset ${\uparrow}\{ \hat{x}, \hat{y}\}$ is closed and finite. Consequently, it forms a finite E-subspace $\boldsymbol{Y}$ of $\X_{n}$. Notice that the restriction $S \coloneqq R \cap (Y \times Y)$ is an E-partition on $\boldsymbol{Y}$. Furthermore, since $\boldsymbol{Y}$ is finite, the E-partition $S$ can be obtained as the kernel of a composition of finitely many $\alpha$ and $\beta$-reductions each of which does not identify any pair of of elements of different color. Consequently, it must be possible to perform on $\boldsymbol{Y}$ at least an $\alpha$ or a $\beta$-reduction $f$ that identifies two distinct elements $z$ and $v$ of the same color. First, suppose that $f$ is a $\beta$-reduction. Then $z$ and $v$ have the same immediate successors. By the definition of the order of $\X_{n}$, if two distinct elements have the same immediate successors, they must be maximal. In particular, $z$ and $v$ are two distinct maximal elements. But this implies that they have different color, a contradiction. Then $f$ must be an $\alpha$-reduction, in which case we can assume, by symmetry, that $v$ is the unique immediate successor of $z$. However, by construction, no element of $\X_{n}$ has precisely one immediate successor, again a contradiction. Hence, we conclude that $\X_{n}$ is $(n+1)$-colorable.
\end{proof}

\begin{Corollary}\label{Cor:n+1-infinite}
Let $n$ be an integer $\geq 2$. The algebra $\X_{n}^{\ast}$ is infinite and $(n+1)$-generated. Consequently, the variety $\VVV(\X_{n}^{\ast})$ is not locally finite.
\end{Corollary}

\begin{proof}
The algebra $\X_{n}^{\ast}$ is $(n+1)$-generated by Theorem \ref{Thm:coloring} and Proposition \ref{Prop:n+1-colorable}. Furthermore, since $\X_{n}$ is infinite, so is $\X_{n}^{\ast}$.
\end{proof}

\section{A combinatorial lemma}

\begin{figure}[h]
\begin{tabular}{ccccccc}
  %\resizebox{0.5\textwidth}{!}{\begin{minipage}{\textwidth}
\\
\begin{tikzpicture}[node distance=8cm and 1cm,auto]
    \tikzstyle{point} = [shape=circle, thick, draw=black, fill=black , scale=0.35]
    \node  (00) at (-5,2) [point,label={\small$y_{0,0}$}] {};
        \node  (01) at (-3,2) [point,label={\small$y_{0,1}$}] {};
        \node  (02) at (-1,2) [point,label={\small$y_{0,2}$}] {};
        \node  (0n) at (5,2) [point,label={\small$y_{0,2^{n+1}-1}$}] {};

 \node (zz1) at (1,2) {$\dots$};
  \node (zz2) at (1,0) {$\dots$};
 \node (zz3) at (1,-2) {$\dots$};
 \node (zz4) at (1,-4) {$\dots$};

     \node  (10) at (-5,0) [point,label={\small$y_{1,0}$}] {};
     \node  (11) at (-3,0) [point,label={\small$y_{1,1}$}] {};
     \node  (12) at (-1,0) [point,label={\small$y_{1,2}$}] {};
     \node  (1n) at (5,0) [point,label={\small$y_{1,2^{n+1}-1}$}] {};
     
     \node  (20) at (-5,-2) [point,label={\small$y_{2,0}$}] {};
     \node  (21) at (-3,-2) [point,label={\small$y_{2,1}$}] {};
     \node  (22) at (-1,-2) [point,label={\small$y_{2,2}$}] {};
     \node  (2n) at (5,-2) [point,label={\small$y_{2,2^{n+1}-1}$}] {};
     
      \node  (30) at (-5,-4) [point,label={\small$y_{3,0}$}] {};
     \node  (31) at (-3,-4) [point,label={\small$y_{3,1}$}] {};
     \node  (32) at (-1,-4) [point,label={\small$y_{3,2}$}] {};
     \node  (3n) at (5,-4) [point,label={\small$y_{3,2^{n+1}-1}$}] {};
    
    \node (dots30a) at (-5,-4.3) {$\vdots$};
    \node (dots30b) at (-5,-4.7) {$\vdots$};
     \node (dots31a) at (-3,-4.3) {$\vdots$};
    \node (dots31b) at (-3,-4.7) {$\vdots$};
    \node (dots32a) at (-1,-4.3) {$\vdots$};
    \node (dots32b) at (-1,-4.7) {$\vdots$};
    \node (dots3na) at (5,-4.3) {$\vdots$};
    \node (dots3nb) at (5,-4.7) {$\vdots$};
    
 %     \node (dots1) at (1,-5) {$\vdots$};
 %  \node (dots2) at (1,-5.4) {$\vdots$};
%\node (dots4) at (1,-6.3)  [point,label={\small$\bot$}] {};

    \draw  (10) -- (01)  ;
  \draw   (10) -- (02);
  \draw   (10) -- (0n) ;
  \draw  (11) -- (00)  ;
  \draw   (11) -- (02);
  \draw   (11) -- (0n) ;
  \draw  (12) -- (00)  ;
  \draw   (12) -- (01);
  \draw   (12) -- (0n) ;
  \draw  (1n) -- (00)  ;
  \draw   (1n) -- (01);
  \draw   (1n) -- (02) ;
  
   \draw  (20) -- (11)  ;
  \draw   (20) -- (12);
  \draw   (20) -- (1n) ;
  \draw  (21) -- (10)  ;
  \draw   (21) -- (12);
  \draw   (21) -- (1n) ;
  \draw  (22) -- (10)  ;
  \draw   (22) -- (11);
  \draw   (22) -- (1n) ;
  \draw  (2n) -- (10)  ;
  \draw   (2n) -- (11);
  \draw   (2n) -- (12) ;

  \draw  (30) -- (21)  ;
  \draw   (30) -- (22);
  \draw   (30) -- (2n) ;
  \draw  (31) -- (20)  ;
  \draw   (31) -- (22);
  \draw   (31) -- (2n) ;
  \draw  (32) -- (20)  ;
  \draw   (32) -- (21);
  \draw   (32) -- (2n) ;
  \draw  (3n) -- (20)  ;
  \draw   (3n) -- (21);
  \draw   (3n) -- (22) ;
      % \draw  (11) -- (00)  (11) -- (02) (11) -- (0n) 

\end{tikzpicture}
\end{tabular}
\caption{The reflexive and transitive closure of $\prec$ on $P_{n}$.}
\label{Fig:Yn}
\end{figure}
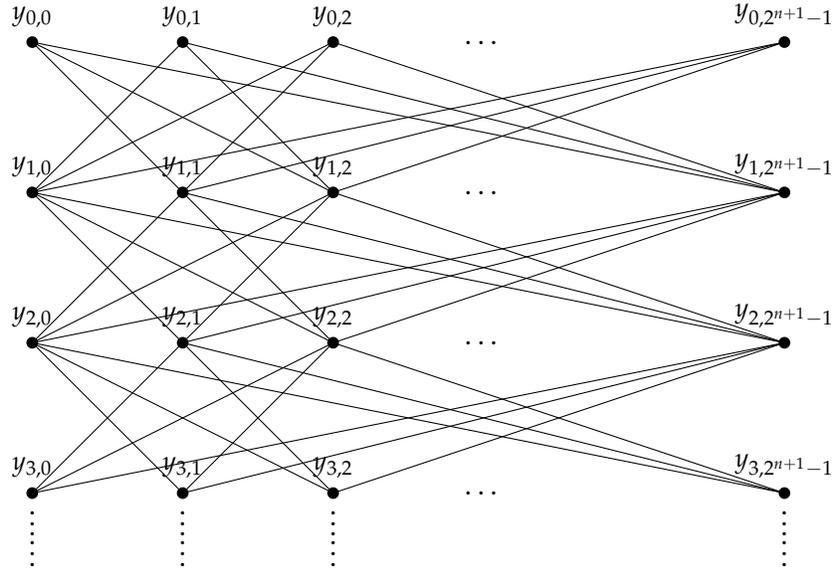

For every $n \in \mathbb{N}$, consider the set
\[
P_{n} \coloneqq \{ y_{m,i} : m,i \in \mathbb{N} \text{ and }i \leq 2^{n+1}-1\}.
\]
We endow $P_{n}$ with a binary relation $\prec$ defined as follows, for every $x, z \in P_{n}$:
\[
x \prec z \Longleftrightarrow x = y_{m,i} \text{ and }z = y_{m-1,j}\text{, for some }m, i, j \in \mathbb{N} \text{ such that }i \ne j.
\]
The poset obtained endowing $P_{n}$ with the reflexive and transitive relation of $\prec$ is depicted in Figure \ref{Fig:Yn}.

\begin{law}
Let $\boldsymbol{Y}_{n}$ be the root compactification of the poset obtained endowing $P_{n}$ with the reflexive and transitive closure of $\prec$. 
\end{law}
\noindent Notice that $\prec$ is the covering relation of the poset $\mathbb{Y}_{n}$ underlying $\boldsymbol{Y}_{n}$. Furthermore, as ${\downarrow}y$ is cofinite for every $y \in Y_{n}$, we can apply of Proposition \ref{Prop:root-compactification}, obtaining that $\boldsymbol{Y}_{n}$ is an Esakia space. 

The aim of this section is to prove the next combinatorial lemma.

\begin{Lemma}\label{Lem:combinatorics}
Let $n \in \mathbb{N}$ and $\boldsymbol{V}$ a finite E-subspace of $\boldsymbol{Y}_{n}$. For every weak $n$-coloring on $\boldsymbol{V}$ there is a finite sequence $\boldsymbol{V}_{0}, \dots, \boldsymbol{V}_{k}$ of Esakia spaces such that:
\benroman
\item\label{item:combinatorics1} $\boldsymbol{V}_{0} = \boldsymbol{V}$ and each $\boldsymbol{V}_{i+1}$ is obtained by applying a $\beta$-reduction $f_{i}$ to $\boldsymbol{V}_{i}$;
\item\label{item:combinatorics2} $\textup{Ker}(f_{k-1} \circ \cdots \circ f_{0})$ does not identify any pair of elements of distinct color of $\boldsymbol{V}$;
\item\label{item:combinatorics3} for every $m \in \mathbb{N}$ such that $y_{m, 0}, \dots, y_{m, 2^{n+1}-1} \in V$, there are $i < j$ such that
\[
\langle y_{m,i}, y_{m,j}\rangle \in \textup{Ker}(f_{k-1} \circ \cdots \circ f_{0}).
\]
\enormal
\end{Lemma}

\begin{proof}
We reason by induction on $n$. For the base case, take $n = 0$. In this case, $\boldsymbol{V}$ is a finite upset of the poset depicted Figure \ref{Fig:Y1}. 
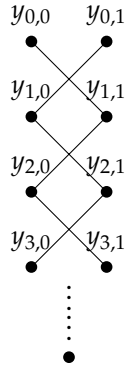
\begin{figure}[h]
\begin{tabular}{ccccccc}
\\
\begin{tikzpicture}
    \tikzstyle{point} = [shape=circle, thick, draw=black, fill=black , scale=0.35]
    \node  (1a) at (-0.5,2) [point,label={\small$y_{0,0}$}] {};
    \node  (1b) at (0.5,2) [point,label={\small$y_{0,1}$}] {};
    \node  (2a) at (-0.5,1) [point,label={\small$y_{1,0}$}] {};
    \node  (2b) at (0.5,1) [point,label={\small$y_{1,1}$}] {};  
  \node  (3a) at (-0.5,0) [point,label={\small$y_{2,0}$}] {};
    \node  (3b) at (0.5,0) [point,label={\small$y_{2,1}$}] {};  
      \node  (4a) at (-0.5,-1) [point,label={\small$y_{3,0}$}] {};
    \node  (4b) at (0.5,-1) [point,label={\small$y_{3,1}$}] {};  
       \node (dots1) at (0,-1.3) {$\vdots$};
       \node (dots2) at (0,-1.7) {$\vdots$};
\node (dots2) at (0,-2.2)  [point] {};
    \draw  (1a) -- (2b) -- (3a) -- (4b)  (1b) -- (2a) -- (3b) -- (4a);

\end{tikzpicture}
\end{tabular}
\caption{The Esakia space $\boldsymbol{Y}_{0}$.}
\label{Fig:Y1}
\end{figure}
Fix a weak $0$-coloring on $\boldsymbol{V}$. In this case, every point of $\boldsymbol{V}$ is colored by the same color, namely the empty sequence of zeros and ones. If there is no nonnegative integer $m \in \mathbb{N}$ such that $y_{m, 0}, y_{m,1} \in V$, then we are done taking $\boldsymbol{V}_{0} = \boldsymbol{V}$ and $k = 0$. Then we consider the case where there is a nonnegative integer $m \in \mathbb{N}$ such that $y_{m, 0}, y_{m,1} \in V$. Let $k$ be the largest such integer (it exists, because $\boldsymbol{V}$ is finite). Define $\boldsymbol{V}_{0} = \boldsymbol{V}$ and let $f_{0} \colon \boldsymbol{V}_{0} \to \boldsymbol{Z}$ be the $\beta$-reduction on $\boldsymbol{V}_{0}$ that identifies $c_{0,0}$ and $c_{0,1}$. We set $\boldsymbol{V}_{1} \coloneqq \boldsymbol{Z}$. Furthermore, let $i$ be a positive integer $< k$ and suppose we already defined $\boldsymbol{V}_{i}$ and $f_{i} \colon \boldsymbol{V}_{i} \to \boldsymbol{V}_{i+1}$. Then let $f_{i+1}$ be the $\beta$-reduction on $\boldsymbol{V}_{i+1}$ that identifies $f_{i} \circ \cdots \circ f_{0}(y_{i+1, 0})$ and $f_{i} \circ \cdots \circ f_{0}(y_{i+1, 1})$ and let $\boldsymbol{V}_{i+2}$ be the codomain of $f_{i+1}$.  Clearly, $\boldsymbol{V}_{0}, \dots, \boldsymbol{V}_{k+1}$ and $f_{0}, \dots, f_{k}$ are ,respectively, sequences of Esakia spaces and $\beta$-reductions satisfying conditions (\ref{item:combinatorics1}) and (\ref{item:combinatorics3}). Finally, condition (\ref{item:combinatorics2}) is satisfied, because all the elements of $\boldsymbol{V}$ have the same color.

For the inductive step, take a positive integer $n$ and suppose that the statement holds for every nonnegative integer $<n$. Let also $\boldsymbol{V}$ be a finite upset of $\boldsymbol{Y}_{n}$ and fix a weak $n$-coloring on $\boldsymbol{V}$. As in the base case, we can assume that there is at least one $m \in \mathbb{N}$ such that $y_{m,0}, \dots, y_{m, 2^{n+1}-1} \in V$. Then set $\boldsymbol{V}_{0} \coloneqq \boldsymbol{V}$. First, we perform a sequence of $\beta$-reductions identifying all maximal elements of $\boldsymbol{V}$, i.e., all elements in $\{y_{0,0}, \dots, y_{0, 2^{n+1}-1}\}$, that have the same color. This gives us a sequence of Esakia spaces $\boldsymbol{V}_{0}, \dots, \boldsymbol{V}_{p_{0}}$ with $\beta$-reductions $f_{0}, \dots, f_{p_{0}-1}$ satisfying conditions (\ref{item:combinatorics1}) and (\ref{item:combinatorics2}), as well as (\ref{item:combinatorics3}) in the restricted case where $m = 0$. 

Now, notice that the elements of depth $\geq 2$ of $\boldsymbol{V}_{p_{0}}$ are the images under $f_{p_{0}-1} \circ \cdots \circ f_{0}$ of the elements of $V$ of the form $y_{i,j}$, for $i \geq 1$. Since $f_{p_{0}-1} \circ \cdots \circ f_{0}$ respects depth (being a composition of $\beta$-reductions), we can assume that the elements of depth $\geq 2$ of $\boldsymbol{V}_{p_{0}}$ are precisely the elements of $V$ of the form $y_{i,j}$ with $i \geq 1$. Furthermore, notice that $\boldsymbol{V}_{p_{0}}$ inherit the weak $n$-coloring of $\boldsymbol{V}$, because the map
\[
f_{p_{0}-1} \circ \cdots \circ f_{0} \colon \boldsymbol{V} \to \boldsymbol{V}_{p_{0}}
\]
does not identify elements of distinct color. 

If it is possible to identify all the elements of depth $2$ of $\boldsymbol{V}_{p_{0}}$ with the same color by means of a series of $\beta$-reduction, we do it and obtain a sequence of finite posets $\boldsymbol{V}_{p_{0}}, \boldsymbol{V}_{p_{0}+1}, \dots, \boldsymbol{V}_{p_{1}}$ and of $\beta$-reductions $f_{p_{0}}, \dots, f_{p_{1}-1}$ such that the sequences of Esakia spaces $\boldsymbol{V}_{0}, \dots, \boldsymbol{V}_{p_{1}}$ and $\beta$-reductions $f_{0}, \dots, f_{p_{1}-1}$ satisfy conditions (\ref{item:combinatorics1}) and (\ref{item:combinatorics2}), as well as (\ref{item:combinatorics1}) in the restricted case where $m \leq 1$. Then, if it is possible, we identify all the elements of depth $3$ of $\boldsymbol{V}_{p_{1}}$ with the same color by means of a series of $\beta$-reduction, we do it and obtain a sequence of finite posets $\boldsymbol{V}_{p_{1}}, \boldsymbol{V}_{p_{1}+1}, \dots, \boldsymbol{V}_{p_{2}}$ and of $\beta$-reductions $f_{p_{1}}, \dots, f_{p_{2}-1}$ such that the sequences of Esakia spaces $\boldsymbol{V}_{0}, \dots, \boldsymbol{V}_{p_{2}}$ and $\beta$-reductions $f_{0}, \dots, f_{p_{2}-1}$ satisfy conditions (\ref{item:combinatorics1}) and (\ref{item:combinatorics2}), as well as (\ref{item:combinatorics1}) in the restricted case where $m \leq 2$. We iterate this process until it is possible and obtain a sequence of Esakia spaces $\boldsymbol{V}_{0}, \dots, \boldsymbol{V}_{p_{k}}$ and $\beta$-reductions $f_{0}, \dots, f_{p_{k}-1}$ that satisfy conditions (\ref{item:combinatorics1}) and (\ref{item:combinatorics2}), as well as (\ref{item:combinatorics1}) in the restricted case where $m \leq k$. If $k$ is the largest nonnegative integer $m$ such that $y_{m,0}, \dots, y_{m,2^{n+1}-1} \in V$, then we are done. Then we consider the case where $y_{k+1,0}, \dots, y_{+1,2^{n+1}-1} \in V$. For the sake of simplicity, as in the case of $\boldsymbol{V}_{p_{0}}$, we can assume that the elements of depth $\geq k+2$ of $\boldsymbol{V}_{p_{k}}$ are precisely the elements of $\boldsymbol{V}$ of the form $y_{i, j}$ with $i \geq k+1$. Consequently, $y_{k+1,0}, \dots, y_{k+1,2^{n+1}-1} \in V_{p_{k}}$.

By assumption, it is not possible to identify all elements in $\{ y_{k+1,0}, \dots, y_{k+1,2^{n+1}-1}\}$ of the same color by applying a series of $\beta$-reductions to $\boldsymbol{V}_{p_{k}}$. This means that there are $i < j$ such that $y_{k+1,i}$ and $y_{k+1,i}$ have the same color, but their immediate successors are different in $\boldsymbol{V}_{p_{k}}$. Since, by construction, $k \geq 0$, we know that the $y_{k+1,i}$ and $y_{k+1,i}$ are not maximal and their immediate successors have depth $k+1$. By the definition of $\boldsymbol{V}_{p_{k}}$ and $\boldsymbol{Y}_{n}$, if the immediate successors of $y_{k+1,i}$ and $y_{k+1,i}$ in $\boldsymbol{V}_{p_{k}}$ are different, then the map $f_{p_{k}-1}\circ \cdots\circ f_{0}$ diverges on $y_{k,i}$ and $y_{k,j}$. In particular, this means that $y_{k,i}$ and $y_{k,j}$ have different color. Therefore, at least one of them is colored by a color different from the constant sequence with value one. By symmetry, we can assume that $y_{k,i}$ is colored by a color $\vec{c} \leq \langle 0, 1, 1, 1 \dots, 1 \rangle$. We shall see that every of depth $\geq k+2$ of $\boldsymbol{V}_{p_{k}}$ is colored by colors $\leq \langle 0, 1, 1, 1 \dots, 1 \rangle$. By definition of $\boldsymbol{Y}_{n}$, it will be enough to prove this for the elements $y_{k+1,0}, \dots, y_{k+1,2^{n+1}-1}$. To this end, notice that all the elements in $\{ y_{k+1,0}, \dots, y_{k+1,2^{n+1}-1}\} \smallsetminus \{ y_{k+1, i}\}$ are below $y_{k,i}$. Thus, since the color of $y_{k,i}$ is $\leq \langle 0, 1, 1, 1 \dots, 1 \rangle$, so are the colors of the elements in $\{ y_{k+1,0}, \dots, y_{k+1,2^{n+1}-1}\} \smallsetminus \{ y_{k+1, i}\}$. Since $y_{k,i}$ and $y_{k,j}$ are different and of the same color, we conclude that also the color of $y_{k,i}$ is $\leq \langle 0, 1, 1, 1 \dots, 1 \rangle$.

Now, let $\vec{c}_{1}, \dots, \vec{c}_{q}$ be the list of distinct colors that are really used to color elements in $\{ y_{k+1,0}, \dots, y_{k+1,2^{n+1}-1}\}$. As we mentioned, $\vec{c}_{i} \leq \langle 0, 1, 1, 1 \dots, 1 \rangle$ for every $i \leq q$. For every $i \leq q$, let $m_{i}$ the the cardinality of the largest set of elements in 
\[
\{ y_{k+1,0}, \dots, y_{k+1,2^{n+1}-1}\} \cap \vec{c}_{i}(\boldsymbol{V}_{p_{k}})
\]
with the same immediate successors in $\boldsymbol{V}_{p_{k}}$. We have two cases:
\benroman
\renewcommand{\labelenumi}{\textbf{\theenumi}}
\renewcommand{\theenumi}{C.\arabic{enumi}}
\item\label{item:greater1} either $m_{1} + \dots + m_{q} \geq 2^{n}$, or
\item\label{item:greater2} $m_{1} + \dots + m_{q} < 2^{n}$.
\eroman

First, suppose that condition (\ref{item:greater1}) holds. In this case, there is a sequence $M_{1}, \dots, M_{t}$ with $t \leq q$ of subsets of $\{ y_{k+1,0}, \dots, y_{k+1,2^{n+1}-1}\}$ such that each $M_{i}$ is a set of elements of color $\vec{c}_{i}$ with the same immediate successors in $\boldsymbol{V}_{p_{k}}$ such that 
\[
\vert M_{1} \vert + \dots + \vert M_{t} \vert = 2^{n}.
\]
Notice that the subposet of $\boldsymbol{Y}_{n}$ with universe
\[
\{ y_{i,j} \in Y_{n} : i \geq k+1 \text{ and }y_{k+1,j} \in M_{1} \cup \dots \cup M_{t} \}
\]
is isomorphic to the poset underlying $\boldsymbol{Y}_{n-1}$. Under this identification, let $\boldsymbol{V}^{-}$ be the finite E-subspace of $\boldsymbol{Y}_{n-1}$ with universe $Y_{n-1} \cap V$. Notice that the set of maximal elements of $\boldsymbol{V}^{-}$ is $M_{1} \cup \dots \cup M_{t}$. Furthermore, the weak $n$-coloring of $\boldsymbol{V}$ restricts to a weak $(n-1)$-coloring of $\boldsymbol{V}^{-}$, because all elements of $\boldsymbol{V}^{-}$ are colored with colors $\leq \langle 0, 1, 1, 1 \dots, 1 \rangle$. 

Therefore, by induction hypothesis, there is a sequence of Esakia spaces $\boldsymbol{V}^{-}_{0}, \dots, \boldsymbol{V}^{-}_{j}$ with $\beta$-reductions $g_{i} \colon \boldsymbol{V}^{-}_{i} \to \boldsymbol{V}^{-}_{i+1}$ for all $j > i \in \mathbb{N}$ such that $\boldsymbol{V}^{-}_{0} = \boldsymbol{V}^{-}$ and:
\benormal
\item\label{item:greater3} $\textup{Ker}(g_{j}\circ \cdots\circ g_{0})$ does not identify any pair of elements of distinct color;
\item\label{item:greater4} for every $m \in \mathbb{N}$ such that $\{ y_{m,i} : y_{k+1,i} \in M_{1} \cup \dots \cup M_{t} \} \subseteq V^{-}$, there are $i < j$ such that $y_{k+1,i}, y_{k+1,j} \in M_{1} \cup \dots \cup M_{t}$ and $\langle y_{m,i}, y_{m,j}\rangle \in \textup{Ker}(g_{j}\circ \cdots\circ g_{0})$. 
\enormal
Now, we can always assume that the $\beta$-reductions $g_{0}, \dots, g_{j}$ are ordered as follows: first we have the $\beta$-reductions $g_{0}, \dots, g_{s_{1}}$ that identify pairs of elements of depth $1$ of $\boldsymbol{V}^{-}$, then those that identify pairs of elements of depth $2$ of $\boldsymbol{V}^{-}$, in symbols $g_{s_{1}+1}, \dots, g_{s_{2}}$, and so on. For each $\beta$-reduction $g_{i}\colon \boldsymbol{V}^{-}_{i} \to \boldsymbol{V}^{-}_{i+1}$, we shall define a $\beta$-reduction
\[
f_{p_{k}+i} \colon \boldsymbol{V}_{p_{k}+i} \to \boldsymbol{V}_{p_{k}+i+1}
\]
as follows. First, $g_{0}$ identifies a pair $\langle x, z \rangle$ of elements of the same color and of depth $1$ in $\boldsymbol{V}^{-}$. Since the set of maximal elements of $\boldsymbol{V}^{-}$ is $M_{1} \cup \dots \cup M_{t}$ and $x$ and $z$ have the same color, we get $x, z \in M_{i}$ for some $i$. By definition, the elements of $M_{i}$ have the same immediate successors in $\boldsymbol{V}_{p_{k}}$, whence the pair $\langle x, z \rangle$ can be identifies by applying a $\beta$-reduction $f_{p_{k}}$ to $\boldsymbol{V}_{p_{k}}$. Essentially the same argument allows to construct the series of $\beta$-reductions $f_{p_{k}}, \dots, f_{p_{k}+s_{1}}$, where $g_{s_{1}+1}$ is the first $\beta$-reduction that identifies a pair of distinct elements $y_{k+2,i}$ and $y_{k+2,j}$ of depth $2$. Then $y_{k+2,i}$ and $y_{k+2,j}$ have the same immediate successors in $\boldsymbol{V}^{-}_{s_{1}+1}$. Since $y_{k+2,i} \prec y_{k+1,j}$ and $y_{k+2,j} \prec y_{k+1,i}$ in $\boldsymbol{V}^{-}_{s_{1}+1}$, we obtain that in $\boldsymbol{V}_{s_{1}+1}^{-}$ the elements $y_{k+1,i}$ and $y_{k+1,j}$ must have been identified, respectively, with some $y_{k+1,h_{i}}$ and $y_{k+1,h_{j}}$ such that $i \ne h_{i}$ and $j \ne h_{j}$. Consequently, $y_{k+1,i}$ and $y_{k+1,j}$ are also identified, respectively, with $y_{k+1,h_{i}}$ and $y_{k+1,h_{j}}$ in $\boldsymbol{V}_{p_{k}+ s_{1}+1}$. Because of the definition of $\boldsymbol{Y}_{n}$, this implies that $y_{k+2,i}$ and $y_{k+2,j}$ have the same immediate successors in $\boldsymbol{V}_{p_{k} + s_{1} + 1}$. Then there exists a $\beta$-reduction $f_{s_{1}+2}$ on $\boldsymbol{V}_{s_{1}+1}$ that identifies the pair $\langle y_{k+2,i}, y_{k+2,j}\rangle$. Essentially the same argument allows to construct the series of $\beta$-reductions $f_{p_{k}+s_{1}+1}, \dots, f_{p_{k}+s_{2}}$. Iterating this process, we obtain a sequence of Esakia spaces $\boldsymbol{V}_{p_{k}+1}, \dots, \boldsymbol{V}_{p_{k}+j+1}$ and of $\beta$-reductions $f_{p_{k}}, \dots, f_{p_{k}+j}$. Clearly, the sequences $\boldsymbol{V}_{0}, \dots, \boldsymbol{V}_{p_{k}+j+1}$ and $f_{0}, \dots, f_{p_{k}+j}$ satisfy condition (\ref{item:combinatorics1}). Together with the definition of the various $f_{i}$, conditions (\ref{item:greater3}) and (\ref{item:greater4}) imply that these sequences satisfy also (\ref{item:combinatorics2}) and (\ref{item:combinatorics3}), as desired. 

Then we consider case where condition (\ref{item:greater2}) holds, i.e., $m_{1} + \dots + m_{q} < 2^{n}$. Since
\[
2m_{1} + \dots + 2m_{q} < 2^{n+1}
\]
and the elements of $\{ y_{k+1,0}, \dots, y_{k+1, 2^{n+1}-1}\}$ are colored with colors among $\vec{c}_{1}, \dots, \vec{c}_{q}$, this implies that there exists $i \leq q$ such that
\[
2m_{i} +1 \leq \vert  \vec{c}_{i}(\boldsymbol{V}_{p_{k}}) \cap \{ y_{k+1,0}, \dots, y_{k+1, 2^{n+1}-1}\} \vert.
\]
Recall that $m_{i}$ is the size of the largest the of elements of $\vec{c}_{i}(\boldsymbol{V}_{p_{k}}) \cap \{ y_{k+1,0}, \dots, y_{k+1, 2^{n+1}-1}\}$ with the same immediate successors. Then, in view of the above display, there are at least three distinct elements $y_{k+1,a_{1}}, y_{k+1,a_{2}}, y_{k+1,a_{3}} \in \vec{c}_{i}(\boldsymbol{V}_{p_{k}}) \cap \{ y_{k+1,0}, \dots, y_{k+1, 2^{n+1}-1}\}$ with different sets of immediate successors. By definition of $\boldsymbol{Y}_{n}$, this means that the elements $y_{k,a_{1}}, y_{k,a_{2}}$, and $y_{k,a_{3}}$ are not identified in $\boldsymbol{V}_{p_{k}}$. By the construction of $\boldsymbol{V}_{p_{k}}$, this guarantees that $y_{k,a_{1}}, y_{k,a_{2}}$, and $y_{k,a_{3}}$ have different colors $d_{1}, d_{2}$, and $d_{3}$. Since $d_{1}, d_{2}, d_{3}$ are distinct, we can assume, without loss of generality, every color $\vec{c} \leq d_{1}, d_{2}, d_{3}$ is $\leq \langle 0, 0, 1, 1 ,1, \dots 1 \rangle$. Notice that every element in $\{y_{k+1,m} \colon m \leq 2^{n+1}-1\} \smallsetminus \{ y_{k+1,a_{1}}, y_{k+1,a_{2}}, y_{k+1,a_{3}} \}$ is below $y_{k,a_{1}}, y_{k,a_{2}}, y_{k,a_{3}}$ and, therefore, of color $\leq \langle 0, 0, 1, 1 ,1, \dots 1 \rangle$. Furthermore, $y_{k+1,a_{1}} \leq y_{k,a_{2}}, y_{k,a_{3}}$ and, therefore, the color of $y_{k+1,a_{1}}$ is $\leq d_{2}, d_{3}$. Similarly, $y_{k+1,a_{2}} \leq y_{k,a_{1}}$, whence the color of $y_{k+1,a_{2}}$ is $\leq d_{1}$. Since $y_{k+1,a_{1}}, y_{k+1,a_{2}}, y_{k+1,a_{3}}$ have the same color, their color is $\leq d_{1}, d_{2}, d_{3}$ and, therefore, $\leq \langle0, 0, 1, 1 ,1, \dots 1 \rangle$. We conclude that every element in $\boldsymbol{V}_{p_{k}}$ of the form $y_{i,j}$ with $i \geq k+1$ is of color $\leq \langle0, 0, 1, 1 ,1, \dots 1 \rangle$. Bearing this in mind, if $m_{1} + \dots + m_{q} \geq 2^{n-1}$, we can conclude the proof by repeating the argument detailed in case (\ref{item:greater1}) with the only different that $n$ should be replaced by $n-1$ in it.

Then we consider the case where $m_{1} + \dots + m_{q} < 2^{n-1}$. There is $i \leq q$ such that
\[
2^{2}m_{i} +1 \leq \vert  \vec{c}_{i}(\boldsymbol{V}_{p_{k}}) \cap \{ y_{k+1,0}, \dots, y_{k+1, 2^{n+1}-1}\} \vert.
\]
As a consequence, there are distinct elements
\[
y_{k+1,a_{1}}, \dots, y_{k+1,a_{5}} \in \vec{c}_{i}(\boldsymbol{V}_{p_{k}}) \cap \{ y_{k+1,0}, \dots, y_{k+1, 2^{n+1}-1}\}
\]
with different sets of immediate successors. By definition of $\boldsymbol{Y}_{n}$, this means that the elements $y_{k,a_{1}}, \dots, y_{k,a_{5}}$ are not identified in $\boldsymbol{V}_{p_{k}}$, whence they have different colors $d_{1}, \dots, d_{5}$. Repeating the argument detailed above paragraph, we obtain that every element in $\boldsymbol{V}_{p_{k}}$ of the form $y_{i,j}$ with $i \geq k+1$ is of color $\leq \langle 0, 0, 0, 1, 1 ,1, \dots 1 \rangle$. Bearing this in mind, if $m_{1} + \dots + m_{q} \geq 2^{n-2}$, we can conclude the proof by repeating the argument detailed in case (\ref{item:greater1}) with the only different that $n$ should be replaced by $n-2$ in it.

As $m_{1} + \dots + m_{q} \geq 1 = 2^{0} = 2^{n-n}$, iterating this argument, eventually we will be able to apply the argument detailed for case (\ref{item:greater1}) and, therefore, conclude the proof.
\end{proof}

\begin{Corollary}\label{Cor:combinatorial-trick}
Let $2 \leq n \in \mathbb{N}$ and let $\boldsymbol{Z}$ be a finite E-subspace of $\X_{n}$ with an E-partition $R$ such that $\boldsymbol{Z}/R$ is $n$-colorable. For every $m \in \mathbb{N}$, if $c_{m, 0}, \dots, c_{m, 2^{n+1}-1} \in Z$, then there are $i < j$ such that $\langle c_{m,i}, c_{m,j}\rangle \in R$.
\end{Corollary}

\begin{proof}
If there is no $m \in \mathbb{N}$ such that $c_{m, 0}, \dots, c_{m, 2^{n+1}-1} \in Z$, the the statement is vacuously true. Then suppose that there is such an integer and let $m$ be the largest one (it exists, because $\boldsymbol{Z}$ is finite). Then let $\boldsymbol{V}$ be the finite E-subspace of $\boldsymbol{Y}_{n}$ with universe
\[
\{ y_{k,i} \in Y_{n} : k \leq 3m \text{ and } i \leq 2^{n+1}-1 \}.
\]
Moreover, let $\delta \colon \boldsymbol{V} \to \boldsymbol{Z}$ be the map defined as follows:
%\begin{displaymath}
%\delta(y_{k,i}) \coloneqq \left\{\begin{array}{@{\,}ll}
%c_{\frac{k}{3},i} & \text{if $k \equiv 0 \mod 3$}\\
%0 & \text{if $a = n$ and $b = 0$}\\
%a & \text{if $b = a - 1$ and $a \geq 3$}\\
%a - 1 & \text{if $b = a - 2$ and $a \geq 3$}\\
%1 & \text{otherwise}\\
%\end{array} \right.
%\end{displaymath}
\benroman
\item $\delta(y_{0,i}) = c_{0,i}$, $\delta(y_{1,i}) = d_{0,i}$, and $\delta(y_{2,i}) = e_{0,i}^{a}$, for every $2^{n+1}-1 \geq i \in \mathbb{N}$;
\item $\delta(y_{3,i}) = c_{1,i}$, $\delta(y_{4,i}) = d_{1,i}$, and $\delta(y_{5,i}) = e_{1,i}^{a}$, for every $2^{n+1}-1 \geq i \in \mathbb{N}$;
\item etc.
\item $\delta(y_{3m,i}) = c_{m,i}$, for every $2^{n+1}-1 \geq i \in \mathbb{N}$.
\eroman
Notice that $\delta$ is a well-defined order embedding. 

Now, recall that $\boldsymbol{Z} / R$ is $n$-colorable and fix an $n$-coloring on it. This $n$-coloring induces a weak $n$-coloring $c$ on $\boldsymbol{Z}$ that colors an element $z \in Z$ by the color of its equivalence class $z / R$. Furthermore, $R$ is the largest E-partition on $\boldsymbol{Z}$ that does not identify any pair of elements of $Z$ colored differently by $c$. In turn, $c$ induces a weak $n$-coloring on $\boldsymbol{V}$ that colors an element $v \in V$ by the color of $\delta(v)$ in $\boldsymbol{Z}$. The fact that this is indeed a weak $n$-coloring on $\boldsymbol{V}$ follows from the fact that $\delta$ is an order embedding.

By Lemma \ref{Lem:combinatorics} there is a finite sequence $\boldsymbol{V}_{0}, \dots, \boldsymbol{V}_{k}$ of Esakia spaces such that:
\benroman
\item $\boldsymbol{V}_{0} = \boldsymbol{V}$ and each $\boldsymbol{V}_{i+1}$ is obtained by applying a $\beta$-reduction $f_{i}$ to $\boldsymbol{V}_{i}$;
\item $\textup{Ker}(f_{k-1} \circ \cdots \circ f_{0})$ does not identify any pair of elements of distinct color of $\boldsymbol{V}$;
\item \label{Eq:combinatorics-3m} for every $3m \geq p \in \mathbb{N}$, there are $i < j$ such that
\[
\langle y_{p,i}, y_{p,j}\rangle \in \textup{Ker}(f_{k-1} \circ \cdots \circ f_{0}).
\]
\enormal
We shall use them to define a sequence $\boldsymbol{Z}_{0}, \dots, \boldsymbol{Z}_{k}$ of Esakia spaces such that:
\benroman
\renewcommand{\labelenumi}{(\theenumi)}
\renewcommand{\theenumi}{\roman{enumi}}
\setcounter{enumi}{3}
\item\label{Eq:Z1-comb} $\boldsymbol{Z}_{0} = \boldsymbol{Z}$ and each $\boldsymbol{Z}_{i+1}$ is obtained by applying a $\beta$-reduction $g_{i}$ to $\boldsymbol{Z}_{i}$;
\item\label{Eq:Z2-comb} $\textup{Ker}(g_{k-1} \circ \cdots \circ g_{0})$ does not identify any pair of elements of distinct color of $\boldsymbol{Z}$;
\item\label{Eq:Z3-comb} for every $m \geq p \in \mathbb{N}$, there are $i < j$ such that
\[
\langle c_{p,i}, c_{p,j}\rangle \in \textup{Ker}(g_{k-1} \circ \cdots \circ g_{0}).
\]
\enormal
To this end, recall that $f_{0} \colon \boldsymbol{V} \to \boldsymbol{V}_{1}$ is a $\beta$-reduction that identifies two elements $x$ and $z$ of the same color. Then $\delta(x)$ and $\delta(z)$ are elements of the same color. Moreover, since $x$ and $z$ have the same successors, they must be maximal, whence so are $\delta(x)$ and $\delta(z)$. Consequently, we can identify $\delta(x)$ and $\delta(z)$ by means of a $\beta$-reduction $g_{0} \colon \boldsymbol{Z} \to \boldsymbol{Z}_{1}$. We repeat this argument, transforming each $f_{i}$ into a $g_{i}$, until we find an $f_{q}$ that identifies two nonmaximal elements $x$ and $z$. Then
\[
x = f_{q-1} \circ \cdots \circ f_{0}(y_{p,i}) \text{ and }z = f_{q-1} \circ \cdots \circ f_{0}(y_{p,j})
\]
for some $p, i, j \in \mathbb{N}$ such that $i \ne j$ and $p \geq 1$. Since $y_{p,i}$ and $y_{p,j}$ are of the same color, so are $\delta(y_{p,i})$ and $\delta(y_{p,j})$. Furthermore, since $x$ and $y$ have the same immediate successors, there are $h_{i} \ne i$ and $h_{j} \ne j$ such that
\begin{align*}
f_{q-1} \circ \cdots \circ f_{0}(y_{p-1,i}) &= f_{q-1} \circ \cdots \circ f_{0}(y_{p-1,h_{i}})\\
f_{q-1} \circ \cdots \circ f_{0}(y_{p-1,j}) &= f_{q-1} \circ \cdots \circ f_{0}(y_{p-1,h_{j}}).
\end{align*}
By definition of $g_{0}, \dots, g_{q-1}$, we get
\begin{align}
g_{q-1} \circ \cdots \circ g_{0}(\delta(y_{p-1,i})) &= g_{q-1} \circ \cdots \circ g_{0}(\delta(y_{p-1,h_{i}}))\label{Eq:the-g-sequence1}\\
g_{q-1} \circ \cdots \circ g_{0}(\delta(y_{p-1,j})) &= g_{q-1} \circ \cdots \circ g_{0}(\delta(y_{p-1,h_{j}})).\label{Eq:the-g-sequence2}
\end{align}
Now, there exists $\hat{p} \in \mathbb{N}$ such that one of the following conditions holds:
\benroman
\renewcommand{\labelenumi}{\textbf{\theenumi}}
\renewcommand{\theenumi}{C.\arabic{enumi}}
\item\label{Eq:caseC1} $\delta(y_{p,i}) = c_{\hat{p},i}$, $\delta(y_{p,j}) = c_{\hat{p},j}$, $\delta(y_{p-1,i}) = e_{\hat{p}-1,i}^{a}$, and $\delta(y_{p-1,j}) = e_{\hat{p}-1,j}^{a}$;
\item $\delta(y_{p,i}) = d_{\hat{p},i}$, $\delta(y_{p,j}) = d_{\hat{p},j}$, $\delta(y_{p-1,i}) = c_{\hat{p},i}$, and $\delta(y_{p-1,j}) = c_{\hat{p},j}$;
\item $\delta(y_{p,i}) = e_{\hat{p},i}^{a}$, $\delta(y_{p,j}) = e_{\hat{p},j}^{a}$, $\delta(y_{p-1,i}) = d_{\hat{p},i}$, and $\delta(y_{p-1,j}) = d_{\hat{p},j}$.
\eroman
We detail only case (\ref{Eq:caseC1}), as the other ones are analogous. By (\ref{Eq:the-g-sequence1}) the elements $e_{\hat{p}-1,i}$ and $e_{\hat{p}-1,h_{i}}$ are identified in $\boldsymbol{Z}_{q}$ (formally, they are identifies by the composition $g_{q-1} \circ \cdots \circ g_{0}$).  Condition (\ref{Eq:the-g-sequence2}) yields the same conclusion for $e_{\hat{p}-1,j}$ and $e_{\hat{p}-1,h_{j}}$. Since $h_{i} \ne i$ and $h_{j} \ne j$, in view of the definition of $\X_{n}$, this implies that the images of $c_{\hat{p},i}$ and $c_{\hat{p},j}$ under $g_{q-1} \circ \cdots \circ g_{0}$ have the same immediate successors in $\boldsymbol{Z}_{q}$. Consequently, we can identify them by means of a $\beta$-reduction $g_{q} \colon \boldsymbol{Z}_{q} \to \boldsymbol{Z}_{q+1}$. This concludes the construction of $g_{q}$.

Repeating this argument, we obtain sequences $\boldsymbol{Z}_{0}, \dots, \boldsymbol{Z}_{k}$ and $g_{0}, \dots, g_{k-1}$ of Esakia spaces and $\beta$-reductions. By construction and (\ref{Eq:combinatorics-3m}) they satisfy conditions (\ref{Eq:Z1-comb}), (\ref{Eq:Z2-comb}), and (\ref{Eq:Z3-comb}), as desired. Finally, take
\[
S \coloneqq \textup{Ker}(g_{k-1} \circ \cdots \circ g_{0}).
\]
As $g_{0}, \dots, g_{k-1}$ is a sequence of $\beta$-reductions, $S$ is an E-partition on $\boldsymbol{Z}$. Furthermore, by (\ref{Eq:Z2-comb}), $S$ does not identify any pair of elements of distinct color. Because $R$ is the largest such E-partition of $\boldsymbol{Z}$, this yields $S \subseteq R$. From (\ref{Eq:Z3-comb}) we conclude that for every $p \in \mathbb{N}$, if $c_{p, 0}, \dots, c_{p, 2^{n+1}-1} \in Z$, then there are $i < j$ such that $\langle c_{p,i}, c_{p,j}\rangle \in R$, as desired.
\end{proof}

\section{The main result}

Our aim is to prove the following:

\begin{Theorem}\label{Thm:main}
For every $n \in \mathbb{N}$ there exists a variety of Heyting algebras whose $n$-generated free algebra is finite, while its $(n+1)$-generated free algebra is infinite.
\end{Theorem}

\noindent The next corollary is an immediate consequence of the theorem:
 
\begin{Corollary}
For every $n \in \mathbb{N}$ there exists a nonlocally finite variety of Heyting algebras whose $n$-generated free algebra is finite.
\end{Corollary}

Theorem \ref{Thm:main} is obvious for the case where $n = 0$, because the variety of all Heyting algebras is nonlocally finite, while its free zero-generated Heyting algebra is finite (being the two-element Boolean algebra). For $n = 1$, it is well known that the variety $\class{KC}$ of Heyting algebras axiomatized by the \textit{weak excluded middle} $\lnot x \lor \lnot \lnot x\thickapprox 1$ is nonlocally finite, while its free one-generated algebra is finite. To prove this, it suffices to observe that every subdirectly irreducible one-generated member of $\class{KC}$ has cardinality $\leq 3$, while the Rieger-Nishimura lattice plus a new bottom element is a two-generated infinite member of $\class{KC}$.

Accordingly, to establish Theorem \ref{Thm:main}, it suffices to exhibit, for every integer $n\geq 2$, a variety of Heyting algebras whose $n$-generated free algebra is finite, while its $(n+1)$-generated free algebra is infinite. In view of Corollary \ref{Cor:n+1-infinite}, the $(n+1)$-generated free algebra of $\VVV(\X_{n}^{\ast})$ is infinite. Therefore, it suffices to prove the following:

\begin{Proposition}\label{Prop:X-n-colorable}
For every integer $n \geq 2$, the $n$-generated free algebra of $\VVV(\X_{n}^{\ast})$ is finite.
\end{Proposition}

\begin{proof}
Since the type of Heyting algebras is finite, it suffices to show that there exists a natural number $m \in \mathbb{N}$ such that the $n$-generated subalgebras of $\X_{n}^{\ast}$ have cardinality $\leq m$. Suppose the contrary, with a view to contradiction. In view of the correspondence between E-partitions and subalgebras and Theorem \ref{Thm:coloring}, there is no natural number $m \in \mathbb{N}$ such that $\vert X_{n} / R\vert \leq m$, for every E-partition $R$ on $\X_{n}$ such that $\X_{n} / R$ is $n$-colorable.

We claim that there are $k \in \mathbb{N}$ and a sequence $\{ \boldsymbol{Z}_{m} : m \in \mathbb{N} \}$ of finite E-subspaces of $\X_{n}$, each with an E-partition $R_{m}$ such that $\boldsymbol{Z}_{m} / R_m$ is $n$-colorable,
\begin{equation}\label{Eq:claim-all-zero}
\vert Z_{m}/ R_{m} \smallsetminus \vec{0}(\boldsymbol{Z}_{m} / R_{m}) \vert \leq k,
\end{equation}
and
\begin{equation}\label{Eq:the-infinite-sequence-Z}
\vert Z_{1}/R_{1} \vert < \vert Z_{2}/R_{2} \vert < \vert Z_{3}/R_{3} \vert < \cdots
\end{equation}

We begin by proving that there is a sequence $\{ \boldsymbol{Y}_{m} : m \in \mathbb{N} \}$ of finite E-subspaces of $\X_{n}$, each with an E-partition $S_{m}$ such that $\boldsymbol{Y}_{m} / S_{m}$ is $n$-colorable and
\begin{equation}\label{Eq:Xn-Y-sequence}
\vert Y_{1}/S_{1} \vert < \vert Y_{2}/S_{2} \vert < \vert Y_{3}/S_{3} \vert < \cdots
\end{equation}
If there is an E-partition $R$ on $\X_{n}$ such that $\X_{n} / R$ is infinite and $n$-colorable, then (because of the definition of $\X_{n}$) there is a sequence
\[
Y_{1} \subsetneq Y_{2} \subsetneq \dots \subsetneq Y_{m} \subsetneq \cdots
\]
of finite upsets of $\X_{n}$ such that
\[
\vert Y_{1}/R \vert < \vert Y_{2}/R \vert < \vert Y_{3}/R \vert < \cdots
\]
Each $Y_{m}$ is finite and, therefore, closed. Thereby, $Y_{m}$ is the universe of an E-subspace $\boldsymbol{Y}_{m}$ of $\X_{n}$. Furthermore, $S_{m} \coloneqq R \cap (Y_{m} \times Y_{m})$ is an E-partition on $\boldsymbol{Y}_{m}$ such that $\boldsymbol{Y}_{m} / S_{m}$ is $n$-colorable and (\ref{Eq:Xn-Y-sequence}) holds. 

Then we consider the case where there is no E-partiton $R$ on $\X_{n}$ such that $\X_{n} / R$ is infinite and $n$-colorable. From the assumption that there is no natural bound on the size of $\X_{n} / R$, provided that it is $n$-colored and $R$ is an E-partition, it follows that there is a sequence of E-partitions $\{ R_{m} : m \in \mathbb{N} \}$ on $\X_{n}$ such that each $\X_{n} / R_{m}$ is finite and $n$-colorable and
\[
\vert X_{n}/R_{1} \vert < \vert X_{n}/R_{2} \vert < \vert X_{n}/R_{3} \vert < \cdots
\]
Given $m \in \mathbb{N}$, let $f \colon \X_{n} \to \X_{n}/ R_{m}$ be the natural Esakia surjection. Since $\X_{n}/ R_{m}$ is finite, we can enumerate its minimal elements as $x_{1} / R_{m}, \dots, x_{k}/ R_{m}$. Furthermore, as the topology of $\X_{n} / R_{m}$ is discrete, the singletons $\{ x_{i}/ R_{m} \}$ are open in $\X_{n}/R_{m}$. Therefore, $f^{-1}[\{x_{i} / R_{m}\}]$ is open in $\X_{n}$. Because of the definition of the topology of $\X_{n}$, this implies that $f^{-1}[\{x_{i}/ R_{m}\}]$ contains an element $z_{i}$ other than $\bot$. Since, in Esakia spaces, principal upsets are closed,
\[
\boldsymbol{Y}_{m} \coloneqq{\uparrow}z_{1} \cup \cdots \cup {\uparrow}z_{k}
\]
is an E-subspace of $\X_{n}$. Moreover, $Y_{m}$ is finite, because each $z_{i}$ is different from $\bot$. Lastly, the restriction of $f$ to $\boldsymbol{Y}_{m}$ is still a surjective Esakia morphism from $\boldsymbol{Y}_{m}$ to $\X_{n} / R_{m}$. Therefore, $S_{m} \coloneqq R_{m} \cap (Y_{m} \times Y_{m})$ is an E-partition of $\boldsymbol{Y}_{m}$ such that $\boldsymbol{Y}_{m} / S_{m}$ is $n$-colorable and (\ref{Eq:Xn-Y-sequence}) holds. 

Recall from Lemma \ref{width of abomination} that the size of antichains in the various $\boldsymbol{Y}_{m}$ is bounded by $2^{n+2}$. Thus, we can apply Lemma \ref{Lem:local-finiteness} to the family $\{ \boldsymbol{Y}_{m} / S_{m} : m \in \mathbb{N} \}$ of $n$-colorable spaces, obtaining a natural number $k$ and a family $\{ \boldsymbol{W}_{m} : m \in \mathbb{N} \}$ of $n$-colorable E-subspaces of spaces in the family $\{ \boldsymbol{Y}_{m}/ S_{m} : m \in \mathbb{N} \}$ such that $\vert W_{m} \smallsetminus \vec{0}(\boldsymbol{W}_{m}) \vert \leq k$, for all $m \in \mathbb{N}$, and
\[
\vert W_{1} \vert < \vert W_{2} \vert < \cdots < \vert W_{m} \vert < \cdots 
\]
For every $m \in \mathbb{N}$, let $Z_{m}$ be the union of the equivalence classes in $W_{m}$. Notice that $Z_{m}$ is a finite upset of $\X_{n}$ and, therefore, the universe of an E-subspace $\boldsymbol{Z}_{m}$ of $\X_{n}$. Furthermore, $R_{m} \coloneqq S_{m} \cap (Z_{m} \times Z_{m})$ is an E-partition on $\boldsymbol{Z}_{m}$ such that $\boldsymbol{Z}_{m} / R_{m}\cong \boldsymbol{W}_{m}$. This concludes the proof of the claim.

Now, recall from the definition of $\langle U_{n}; \prec\rangle$ that $T_{n} = \{ s_{0}, \dots, s_{t} \}$. To conclude the proof, it suffices to show that for every $m \in \mathbb{N}$,
\begin{equation}\label{Eq:contradiction}
\vert Z_{m}/ R_{m} \vert \leq k+1 + (3+t)(2^{n+3}+2).
\end{equation}
This is because the above display is in contradiction with (\ref{Eq:the-infinite-sequence-Z}). To this end, fix an arbitrary $m \in \mathbb{N}$ and an $n$-coloring on $\boldsymbol{Z}_{m} / R_{m}$ which satisfies (\ref{Eq:claim-all-zero}). For every $p \in \mathbb{N}$, let $V_{p}$ be the subset of $U_{n}$ consisting of all elements of the form $a_{p}, b_{p}, c_{p, i}, d_{p, i}, e_{p, i}^{a}, e_{p, i}^{b}$, where $i$ ranges over all natural numbers. Clearly,
\[
U_{n} = \bigcup_{p \in \mathbb{N}}V_{p}.
\]
 Now, if there is no $q\in \mathbb{N}$ such that
\begin{equation}\label{Eq:distinct-zero-Z}
\Big((Z_{m} \cap V_{q}) / R_{m}\Big) \cap \vec{0}(\boldsymbol{Z}_{m} / R_{m}) \ne \emptyset,
\end{equation}
then all the elements of $\boldsymbol{Z}_{m} / R_{m}$ are colored with colors other than $\vec{0}$. Hence, by (\ref{Eq:claim-all-zero}), $\vert Z_{m} / R_{m} \vert \leq k$ and we are done.

Then, suppose that (\ref{Eq:distinct-zero-Z}) holds for some $q \in \mathbb{N}$. We can further assume that $q$ is the least natural number validating it. Bearing this in mind, it is clear that
\[
(\bigcup_{q-1 \geq p \in \mathbb{N}}V_{p} \cap Z_{m}) / R_{m} \cap \vec{0}(\boldsymbol{Z}_m/R_m) = \emptyset,
\]
hence it now follows from (\ref{Eq:claim-all-zero}) that 
\begin{equation}\label{Eq:below-k-2}
\vert (\bigcup_{q-1 \geq p \in \mathbb{N}}V_{p} \cap Z_{m}) / R_{m} \vert \leq k.
\end{equation}

Since each $V_{i}$ has cardinality $2^{n+3}+2$, from (\ref{Eq:below-k-2}) it follows
\begin{align*}
\vert (\bigcup_{q+t+2\geq p \in \mathbb{N}}V_{p} \cap Z_{m})/ R_{m} \vert & \leq \vert (\bigcup_{q-1 \geq p \in \mathbb{N}}V_{p} \cap Z_{m}) / R_{m} \vert  + \vert V_{q} \vert + \dots + \vert V_{q+t+2} \vert\\
&\leq k + \vert V_{q} \vert + \dots + \vert V_{q+t+2} \vert\\
& = k+ (3+t)(2^{n+3}+2).
\end{align*}
Therefore, in order to prove (\ref{Eq:contradiction}), it suffices to show that
\begin{equation}\label{Eq:display8}
\vert (\bigcup_{q+t+3 \leq p \in \mathbb{N}} V_{p}\cap Z_{m}) / R_{m} \vert \leq 1.
\end{equation}

To this end, recall that at least one element of $(Z_{m} \cap V_{q}) / R_{m}$ is colored by $\vec{0}$. By the construction of $\X_{n}$, every element of the form $d_{q+1, i}$ is below every element of $V_{q}$. Consequently, every element of the form $d_{q+1, i} / R_{m}$ has color $\vec{0}$ in $\boldsymbol{Z}_{m}/ R_{m}$. Furthermore, every element in $V_{p}$, for $p \geq q+2$, is below every element of the form $d_{q+1, i}$ in $\X_{n}$, whence
\begin{equation}\label{Eq:display7}
(\bigcup_{q+2\leq p \in \mathbb{N}}V_{p} \cap Z_{m}) / R_{m} \subseteq \vec{0}(\boldsymbol{Z}_{m} / R_{m}).
\end{equation}

We have two cases: either $V_{q+t+2} \subseteq Z_{m}$ or $V_{q+t+2} \nsubseteq Z_{m}$.  First suppose that $V_{q+t+2} \nsubseteq Z_{m}$. Since $Z_{m}$ is an upset of $\X_{n}$ and
\[
V_{q+t+2} \subseteq {\uparrow} \{ c_{q+t+3, i}, c_{q+t+3, j}\}, \text{ for every } i < j,
\]
we conclude that $Z_{m}$ contains at most one element of the form $c_{q+t+3, i}$. As every element in $V_{q+t+3}$ is either of the form $c_{q+t+3,i}$ or is below at least two such elements, we conclude that $Z_{m} \cap V_{q+t+3}$ is either empty or a singleton of the form $\{ c_{q+t+3, i} \}$. Together with the fact that $Z_{m}$ is an upset of $\X_{n}$ and the definition of $\X_{n}$, this implies
\[
\vert \bigcup_{q+t+3 \leq p \in \mathbb{N}} V_{p}\cap Z_{m} \vert \leq 1
\]
and, therefore, (\ref{Eq:display8}).

Then we consider the case where $V_{q+t+2} \subseteq Z_{m}$. Since $Z_{m}$ is an upset of $\X_{n}$, this implies $V_{p} \subseteq Z_{m}$ for every $q+t+2 \geq p \in \mathbb{N}$. In particular, $V_{q+2} \subseteq Z_{m}$. Since $\boldsymbol{Z}_{m} / R_{m}$ is $n$-colorable, we can apply Corollary \ref{Cor:combinatorial-trick} obtaining that there are $i < j$ such that $\langle c_{q+2,i}, c_{q+2,j}\rangle \in R_{m}$. We shall prove that
\begin{equation}\label{Eq:the-levels-ZR}
\langle c_{p,i}, c_{p,j}\rangle \in R_{m}, \text{ for every }q+2\leq p \leq q+ t+2.
\end{equation}
The base case, where $p = q+2$, was already established. Then suppose the above display holds for some $q+2\leq p < q+2+t$. Recall from (\ref{Eq:display7}) that the element $c_{p,i} / R_{m} = c_{p,j} / R_{m}$ is colored by $\vec{0}$ in $\boldsymbol{Z}_{m} / R_{m}$. Since $\langle c_{p,i}, c_{p,j}\rangle \in R_{m}$, the elements $d_{p,i} / R_{m}$ and $d_{p,j} / R_{m}$ have the same immediate successors in $\boldsymbol{Z}_{m} / R_{m}$. Therefore, since $\boldsymbol{Z}_{m} / R_{m}$ is finite, there is a $\beta$-reduction that identifies only $d_{p,i} / R_{m}$ and $d_{p,j} / R_{m}$. As these elements are colored by $\vec{0}$, by (\ref{Eq:display7}), this $\beta$-reduction does not identify any pair of elements of distinct color. Together with the fact that $\boldsymbol{Z}_{m} / R_{m}$ is $n$-colorable, this implies that $\langle d_{p,i}, d_{p,j}\rangle \in R_{m}$. Repeating this argument, one shows that also $\langle e_{p,i}^{a}, e_{p,j}^{a}\rangle \in R_{m}$ and, finally, $\langle c_{p+1,i}, c_{p+1,j}\rangle \in R_{m}$, thus establishing (\ref{Eq:the-levels-ZR}).

Now, recall that $T_{n} = \{ s_{0}, \dots, s_{t}\}$. Choose an element $s_{g} \in T_{n}$ of the form $\langle h, i, j \rangle$ (this is possible, because $i \ne j$). Then let $v$ be the unique integer in the interval $[q+2, q+2+ t]$ such that $v \equiv g \mod t+1$. Then
\benroman
\item the immediate successors of $a_{v}$ in $\X_{n}$ are $c_{v,h}$ and $c_{v,i}$;
\item the immediate successors of $b_{v}$ in $\X_{n}$ are $c_{v,h}$ and $c_{v,j}$.
\eroman
Since $\langle c_{v,i}, c_{v,j}\rangle \in R_{m}$, the elements $a_{v}/R_{m}$ and $b_{v}/ R_{m}$ have the same immediate successors in $\boldsymbol{Z}_{m} / R_{m}$. Therefore, since $\boldsymbol{Z}_{m} / R_{m}$ is finite, there is a $\beta$-reduction that identifies only $a_{v} / R_{m}$ and $b_{v} / R_{m}$. As these elements are colored by $\vec{0}$, by (\ref{Eq:display7}), this $\beta$-reduction does not identify any pair of elements of distinct color. Together with the fact that $\boldsymbol{Z}_{m} / R_{m}$ is $n$-colorable, this implies that $\langle a_{v}, b_{v}\rangle \in R_{m}$. Because of this, the elements of the form $e_{v,p}^{a} / R_{m}$ and $e_{v,p}^{b} / R_{m}$ have the same immediate successors in $\boldsymbol{Z}_{m} / R_{m}$. Again, there is a $\beta$-reduction that identifies only $e_{v,p}^{a} / R_{m}$ and $e_{v, p}^{b} / R_{m}$. As these elements are colored by $\vec{0}$, by (\ref{Eq:display7}), this $\beta$-reduction does not identify any pair of elements of distinct color. Together with the fact that $\boldsymbol{Z}_{m} / R_{m}$ is $n$-colorable, this implies that $\langle e_{v,p}^{a}, e_{v,p}^{b}\rangle \in R_{m}$. Thus,
\[
\langle e_{v,p}^{a}, e_{v,p}^{b}\rangle \in R_{m}\text{, for all }2^{n+1}-1\geq p \in \mathbb{N}
\]
Since every element of the form $c_{v+1, u}$ is below every element of the form $e_{v,p}^{b}$ in $\X_{n}$, the above display guarantees that every element of the form $c_{v+1, u} / R_{m}$ is below every element of the form $e_{v,p}^{a} / R_{m} = e_{v,p}^{b} / R_{m}$. As a consequence, the elements of the following sequence have the same immediate successors in $\boldsymbol{Z}_{m} / R_{m}$:
\[
c_{v+1,0} / R_{m}, \dots, c_{v+1, 2^{n+1}-1}/ R_{m}.
\]
Notice that, however, that the above notation contains a minor abuse of notation, $Z_{m}$ need not include the whole $\{ c_{v+1,0}, \dots, c_{v+1, 2^{n+1}-1} \}$.

Repeating once more the argument involving $\beta$-reductions, we obtain that 
\[
c_{v+1,0} / R_{m} =  \cdots = c_{v+1, 2^{n+1}-1}/ R_{m}.
\]
We fix the notation $x \coloneqq c_{v+1,0} / R_{m}$. In view of the above display, either $a_{v+1}/ R_{m} = x$, or $x$ is the only immediate successor of $a_{v+1} / R_{m}$. In the latter case, we could perform a proper $\alpha$-reduction on $\boldsymbol{Z}_{m} / R_{m}$, identifying the elements $a_{v+1} / R_{m}$ and $x$ that, moreover, have the same color by (\ref{Eq:display7}). But this would contradict the fact that $\boldsymbol{Z}_{m} / R_{m}$ is $n$-colorable. Then we conclude that $a_{v+1}/ R_{m} = x$. A similar argument shows that $b_{v+1}/ R_{m} = x$ and $d_{v+1, u} / R_{m} =  x$ for all $2^{n+1}-1 \geq u \in \mathbb{N}$. Since
\[
a_{v+1}/ R_{m}  = b_{v+1}/ R_{m} = d_{v+1, u} / R_{m} =  x\text{, for all }2^{n+1}-1 \geq u \in \mathbb{N},
\]
either every element of the form $e_{v+1, u}^{a} / R_{m}$ is equal to $x$, or $x$ is the only immediate successor of $e_{v+1, u}^{a} / R_{m}$. Repeating the same argument, involving $\alpha$-reductions, we obtain $e_{v+1, u}^{a} / R_{m} = x$. Similarly, every element of the form $e_{v+1, u}^{b} / R_{m}$ is equal to $x$. Thus, we conclude that $(V_{v+1} \cap Z_{m})/ R_{m} \subseteq \{ x \}$. Iterating this argument, we obtain 
\[
(\bigcup_{ v+1 \leq p \in \mathbb{N}}V_{p} \cap Z_{m}) / R_{m} \subseteq \{ x \}
\]
and, therefore,
\[
\vert (\bigcup_{ v+1 \leq p \in \mathbb{N}}V_{p} \cap Z_{m}) / R_{m} \vert \leq 1.
\]
Since $v \leq q+t+2$,  this implies (\ref{Eq:display8}), thus concluding the proof.
\end{proof}

\vspace{.5cm}
\paragraph{\bfseries Acknowledgements.}
%We gratefully thank the anonymous referee for the careful reading of the manuscript.
The first author was supported by the research grant PREDOCS-UB 2021 funded by the University of Barcelona. The second author was supported by the \emph{Beatriz Galindo} grant BEAGAL\-$18$/$00040$ funded by the Ministry of Science and Innovation of Spain. All the authors were supported by the MSCA-RISE-Marie Skłodowska-Curie Research and Innovation Staff Exchange (RISE) project MOSAIC $101007627$ funded by Horizon $2020$ of the European Union.

\bibliographystyle{plain}

\end{document}